\documentclass[preprint,12pt]{amsart}
\usepackage{amssymb,amsfonts,amsmath,amsthm,verbatim}
\usepackage{xcolor}
\usepackage{tikz}
\usepackage{graphicx}
\usepackage[left=2.5cm, right=2.5cm, top=3.5cm, bottom=3.5cm]{geometry}
\usepackage{graphicx}
\usepackage{hyperref}
\usepackage{mathdots}
\usepackage{diagbox}
\usepackage[normalem]{ulem}
\usepackage{array,tabularx} 
\usepackage{booktabs}
\usepackage{multirow}
\usepackage{float}
\newtheorem{theorem}{Theorem}[section]
\newtheorem{theo}[theorem]{Theorem}
\newtheorem{prop}[theorem]{Proposition}
\newtheorem{lem}[theorem]{Lemma}

\newtheorem{defn}[theorem]{Definition}

\newtheorem{exam}[theorem]{Example}

\makeatletter \@addtoreset{equation}{section}

\title{Ehrhart Series and Transfer Matrices for Period-Two Polytopes with Consecutive Inequalities}

\author{Yangbiao Zhou}
\thanks{*Corresponding author}
\address{School of Mathematics and Statistics, Hainan University, Haikou 570228,
P.R. China}

\email{zhouyangbiao520@gmail.com}

\date{\today}
\begin{document}

\begin{abstract}
We study a family of period-two polytopes that can be naturally realized as lattice polytopes. In contrast to polytopes defined by a uniform system of inequalities, this family is characterized by consecutive coordinate inequalities whose bounds alternate between $s+1$ and $s$. The lattice points of these polytopes admit an explicit transfer-matrix formulation, which provides an effective framework for their enumeration and analysis. Furthermore, by considering fixed dilations, we investigate the associated dimensional generating functions and derive rational forms for the generating functions of the lattice points in the $m$-th dilates of the polytopes.
\end{abstract}

\maketitle

\vspace{-5mm}

\noindent
\begin{small}
 \emph{Mathematic subject classification}: 52B05, 05A15,05A19,11D75
\end{small}

\noindent
\begin{small}
\end{small}

\section{introduction}
For a polytope, there are well-established theoretical frameworks for studying the enumeration of its lattice points. In particular, for lattice polytopes, Ehrhart's theorem \cite{Ehrhart1,Ehrhart2} provides a fundamental and effective tool for analyzing their lattice-point enumerators. Jiang-Wen-Zhong \cite{Jiang} studied an alternating adjacent-sum polytope defined by consecutive-coordinate inequalities whose bounds alternate between $s$ and $s+1$. In this paper, we investigate how the corresponding results change under duality. To this end, we introduce another period-two family of adjacent-sum lattice polytopes and study its lattice-point enumeration and associated generating functions.
\[
\mathcal{P}_d^{(s)}=\{(x_a,\dots,x_d)\in \mathbb{R}_{\geq0}^d:x_i+x_{i+1}\leq s+\delta_i,1\leq i\leq d-1\}.
\]
where $\delta_i=1$ if $i$ is odd;$\delta_i=0$ if $i$ is even.Consecutive-coordinate inequalities of this type have been studied from several perspectives. Bona-Ju \cite{Bona} enumerated the solutions to homogeneous linear inequalities motivated by magic squares, while Bona-Ju-Yoshida \cite{Bona-Ju} extended their approach to related problems involving weighted graphs. Xin-Zhong \cite{Zhong} further obtained closed-form identities for the corresponding Kekule numbers. These developments suggest that transfer-matrix methods provide a natural and effective framework for studying enumeration problems governed by consecutive-coordinate constraints. In particular, Engström-Kohl \cite{transfer} explored several classes of combinatorial problems arising in the context of Ehrhart theory that admit a transfer-matrix formulation.Clearly,this polytope is an integer convex $d$-polytope.And when $m=1$,the counts of lattice points of it denoted by $N_d(s)$,i.e.$N_d(s)=\sharp(\mathcal{P}_d^{(s)}\cap \mathbb{Z}^d)$.For instance,when $s=1$ and $d=4$,the constraints are $x_i+x_{i+1}\leq1+\delta_i$ for $i=1,2,3$.Then
\[
N_4(1)=\sum_{x_2=0}^1\sum_{x_3=0}^{1-x_2}(3-x_2)(3-x_3)=21.
\]
While for comparison,the uniform polytope with all bounds equal to 1 possess $\sharp\{\mathbf{x}\in \mathbb{Z}_{\geq 0}^4:x_1+x_{i+1}\leq1,i=1,2,3\}=8$ lattice points.When $d=5$,we find that the constraint $x_3+x_4\leq2$ must be satisfied only if $x_2+x_3\leq1$ and $x_4+x_5\leq1$.Hence the polytope $\mathcal{P}_5^{(1)}$ can be formed by
\[
x_1+x_2\leq2,x_2+x_3\leq1,x_4+x_5\leq1.
\]
Naturely this polytope have $N_5(1)=(\sum_{x_2=0}^1(3-x_2)(2-x_2))(\sum_{x_4=0}^1(2-x_4))=24$ lattice points,but the uniform have only $13$ lattice points.For this conditions,we discuss the value of $d$ as $odd$ or $even$ respectly,given initial value $N_1(s)=s+2$ and defines
\[
\Phi_s(y)=N_1(s)+\sum_{m\geq1}N_{2m+1}(s)y^m,\Psi_s(y)=\sum_{m\geq0}N_{2m+2}(s)y^m.
\]
The main results as follows:
\begin{theo}\label{pqF}
The generating function of $N_d(s)$ in $d$ denoted by $F(x)=\sum_{d\geq2}N_d(s)x^d$,then
\[
F(x)=x^2\Psi_s(x^2)+x(\Phi_s(x^2)-N_1(s)),
\]
where $N_1(s)=s+2$ and $\Psi_s(y),\Phi_s(y)$ have closed-expression for $s\geq1$ as follows
\[
\Phi_s(y)=\frac{P_{s+1}(y)}{Q_{s+2}(y)},\Psi_s(y)=\frac{(1+y)^{s+2}-Q_{s+2}(y)}{yQ_{s+2}(y)}.
\]
\end{theo}
Theorem \ref{pqF} induces the following recurrence between odd and even parts series,
\begin{theo}\label{pqR}
For integer $s\geq1$,the first-order recursive relation:
\begin{gather*}
\Phi_s(y)=\frac{1+\Phi_{s-1}(y)}{1-y\Phi_{s-1}(y)},\Phi_0(y)=\frac{2}{1-y},\\
1+y\Phi_s(y)=\frac{1+y\Psi_s(y)}{1+y\Psi_{s-1}(y)},\Psi_0(y)=\frac{y+3}{1-y}.
\end{gather*}
While the second-order recursive relation for $s\geq2$:
\begin{gather*}
\Phi_s(y)=\frac{2+(1-y)\Phi_{s-2}(y)}{1-y-2y\Phi_{s-2}(y)},\Phi_0(y)=\frac{2}{1-y},\Phi_1(y)=\frac{3-y}{1-3y},\\
\frac{1+y}{1+y\Psi_s(y)}=\frac{2}{1+y\Psi_{s-1}(y)}-\frac{1}{1+y\Psi_{s-2}(y)},\Psi_0(y)=\frac{y+3}{1-y},\Psi_1(y)=\frac{6+3y+y^2}{1-3y}.
\end{gather*}
\end{theo}

For better comparison,addding a extra constraint $x_d+x_1\leq s+\delta_d$ and defining $N_d^{cyc}(s)$the counts of lattice points of the new polytope,all indices are interpreted modulo $d$,\[N_d^{cyc}(s)=\{(x_1,\dots,x_d)\in \mathbb{Z}_{\geq0}^d:x_i+x_{i+1}\leq s+\delta_i,1\leq i\leq d\}.\]
Samely we use the two generating functions in odd- and even-conditions,
\[
\Phi_s^{cyc}(y)=\sum_{m\geq1}N_{2m+1}^{cyc}(s)y^m,\Psi_s^{cyc}(y)=\sum_{m\geq1}N_{2m}^{cyc}(s)y^m
\]
\begin{theo}[cyclic case]\label{pqC}
The generating function be denoted by $F^{cyc}(x)$,and
\[
F^{cyc}(x)=\sum_{d\geq2}N_d^{cyc}(s)x^d=\Psi_s^{cyc}(x^2)+x\Phi_s^{cyc}(x^2),
\]
where $\Phi_s^{cyc}(y)$ and $\Psi_s^{cyc}(s)$ have the following forms for $s\geq0$,
\[
\Phi_s^{cyc}(y)=\frac{\sum_{i=0}^{s+}U_{s+1-i,i}-\sum_{j=0}^sU_{j,s-j}}{Q_{s+2}(y)}-\left[\frac{s+3}{2}\right],
\Psi_s^{cyc}(y)=\frac{y(s+2)P_{s}(y)}{2Q_{s+2}(y)}.
\]
\end{theo}
About the Ehrhart theory,we use a generating function in demension $d$.Firstly,defining the enumerations of the lattice points of the delites of polytope $\mathcal{P}_d^{(s)}$ as $L_d^{(s)}(m)$ for $d\geq2$,thus
\[
L_d^{(s)}(m)=\sharp(m\mathcal{P}_d^{(s)}\cap \mathbb{Z}^d),
\]
And let us assuming $L_0^{(s)}(m)=1,L_1^{(s)}(m)=(s+1)m+1$ and the generating function of the Ehrhart polynomial
\[
G_m^{(s)}(y)=\sum_{d\geq0}L_d^{(s)}(m)y^d.
\]
\begin{theo}\label{Gmy}
For variable $s,m\geq1$,the generating function of the lattice points of the polytope,
\[
G_m^{(s)}(y)=1+y^2U_m^{(s)}(y^2)+yV_m^{(s)}(y^2).
\]
where
\begin{gather*}
U_m^{s}(z)=\frac{\bar{\beta}_0^{(s,m)}(\bar{\beta}_0^{(s,m)}+1)}{2}+\frac{z((\bar{\beta}_0^{(s,m)}+1-q_0(z))p_0(z)\eta_1^{(s,m)}-q_0(z)\xi_1^{(s,m)}(1-zp_0(z)))}{1-zp_0(z)},\\
V_m^{(s)}(z)=\bar{\beta}_0^{(s,m)}+\frac{z((\bar{\beta}_0^{(s,m)}+1)p_0(z)-q_0(z))}{1-zp_0(z)}.
\end{gather*}
where $\eta_j^{(s,m)}=\frac{j((2s+3)m+3)}{2j+1}$ and $\xi_j^{(s,m)}=\frac{((2j+1)s+3j+1)m+3j+1}{2j+2}$.

And for two speacial cases:\\
Case1:fixed $m=1$ and for every $s\geq1$,
\[
C_1^{(s)}(y)=1+(s+2)y+F(y).
\]
Case2:fixed $s=1$ and for every $m\geq1$,
\[
G_m^{(1)}(y)=\frac{1+(2m+1)y+3\binom{m+1}{2}y^2-\binom{m+2}{3}y^3+\frac{9m+3}{4}\binom{m+2}{3}y^4-\frac{1}{3}\binom{m+2}{2}\binom{m+3}{4}y^6}{1-\binom{m+2}{2}y^2}.
\]
\end{theo}

\section{transfer matrices and dimension generating function}

In this section, we fix \(s\geq 1\) and investigate the dependence of the generating functions on the dimension \(d\), treating the cases of odd and even \(d\) separately. For comparison with the uniform model, we also consider the variant obtained by imposing the additional cyclic constraint
\[
x_d+x_1\leq s+\delta_d.
\]
We denote its generating function, viewed as a function of \(d\), by \(F^{\mathrm{cyc}}(x)\).
.
\subsection{transfer matrix of the inequality system}

Given an inequality system $x_k+x_{k+1}\leq s+\delta_k$ for $1\leq k\leq d-1$,where $\delta_k=1$ if $k$ is odd and otherwise $\delta_k=0$.For index $0\leq i,l\leq s+1$,let $a_{ij}(d,s)$ denote the counts of nonnegative integer vectors $(x_1,\dots,x_d)$ satisfying this system with $x_1=i$ and $x_d=j$.And write a matrix $A_{d,s}:=(a_{ij}(d,s))_{0\leq i,j\leq s+1}$.Thus the number of the nonnegative integer solution of the system is
\begin{align}
N_d(s)=\sum_{i=0}^{s+1}\sum_{j=0}^{s+1}a_{ij}(d,s).\label{Nds}
\end{align}

Define two transfer matrices $B_s$ and $\bar{B}_s$ of size $(s+2)\times(s+2)$,but have difference with traditional-matrix index $\{0,1,\dots,s+1\}$,by
\[
(B_s)_{ij}=\begin{cases}
1,&\text{ }i+j\leq s,\\
0,&\text{otherwise}.
\end{cases},
(\bar{B}_s)_{ij}=\begin{cases}
1,&\text{ }i+j\leq s+1,\\
0,&\text{otherwise}.
\end{cases}.
\]
For instance let $s=4$,
\[
B_4=\begin{pmatrix}
1 &1 &1 &1 &1 &0\\
1 & 1 &1 &1&0 &0\\
1& 1& 1&0 &0 &0\\
1&1 &0 &0 &0 &0\\
1 &0 &0 &0 &0 &0\\
0 &0 &0 &0 &0 &0
\end{pmatrix},
\bar{B}_4=\begin{pmatrix}
1&1 &1 &1 &1 &1 \\
1 & 1 &1 &1&1 &0\\
1& 1& 1&1 &0 &0\\
1&1 &1 &0 &0 &0\\
1 &1 &0 &0 &0 &0\\
1 &0 &0 &0 &0&0
\end{pmatrix},
\]
Here we can see that $B_s$ records transitions across the enen-indexed constraint $x_{2r}+x_{2r+1}\leq s$,and $\bar{B}_s$ records the odd-indexed constraint.Let the double-step transfer matrix is
\[
D_s=\bar{B}_sB_s.
\]

Set $u_s\in \mathbb{R}^{s+2}$ is all-one volume vector and write $N_d(s)=u_s^TA_{d,s}u_s$.Thus the closed-expression of odd and even,
\begin{prop}\label{N2msN2m+1s}
For any integer $m\geq1$,
\[
N_{2m}(s)=u_s^TD_s^{m-1}\bar{B}_su_s,N_{2m+1}(s)=u_s^TD_s^mu_s.
\]
\end{prop}

\begin{proof}
By the two tranfer matrices,thus for $m\geq1$ and $d=2m$ is a even,the system have $2m-1$ constraints,samely for $d=2m+1$ have $2m$ constraints,then
\[
A_{2m,s}=\bar{B}_sB_s\cdots\bar{B}_s=D_s^{m-1}\bar{B}_s,
A_{2m+1,s}=\bar{B}_sB_s\cdots\bar{B}_sB_s=D_s^{m}.
\]
By equation \ref{Nds},$N_d(s)$ is the sum of all entries of matrix $A_{d,s}$,it can be written $N_d(s)=u_s^TA_{d,s}u_s$.Then substituting $A_{d,s}$ into this formula,
\[
N_{2m}(s)=u_s^TA_{2m,s}u_s=u_s^T(D_s)^{m-1}\bar{B}_su_s,
N_{2m+1}(s)=u_s^TA_{2m+1,s}u_s=u_s^T(D_s)^mu_s.
\]
This proof the proposition.
\end{proof}
\begin{lem}\label{phipsi}
For $s\geq1$
\[
\Phi_s(y)=\frac{u_s^T(I-yD_s)^*u_s}{det(I-yD_s)},\Psi_s(y)=\frac{u_s^T(I-yD_s)^*\bar{B}_su_s}{det(I-yD_s)}.
\]
where $I$ is $n$-dimensional iidentity matrix and $M^*$ denote adjugate matrix of $M$.
\end{lem}
\begin{proof}
By proposition \ref{N2msN2m+1s} and definition of $\Phi_s(y)$ and $\Psi_s(y)$,
\begin{gather*}
\Phi_s(y)=N_0(y)+\sum_{m\geq1}N_{2m+1}(s)y^m=\sum_{i\geq0}u_s^T(D_s)^mu_sy^m=u_s^T(I-yD_s)^{-1}u_s,\\
\Psi_s(y)=\sum_{m\geq0}N_{2m+2}(s)y^m=\sum_{m\geq0}u_s^T(D_s)^{m}\bar{B}_su_sy^m=u_s^T(I-yD_s)^{-1}\bar{B}_su_s.
\end{gather*}
This adjugate form follows from
\[
(I-yD_s)^{-1}=\frac{(I-yD_s)^*}{det(I-yD_s)}.
\]
Thus substitute it into above formulas such that prove this lemma.
\end{proof}
\subsection{the closed-expression of $\Phi_s(y),\Psi_s(y)$ with the two polynomials $P_{n}(y),Q_{n}(y)$.}In this subsection, we establish that the class polynomials satisfy a nontrivial recurrence in \(n\). We further investigate whether the two adjugate numerators,
\[
u_s^{T}(I-yD_s)^{*}u_s
\qquad\text{and}\qquad
u_s^{T}(I-yD_s)^{*}\bar{B}_s u_s,
\]
satisfy the same recurrence. We also examine the denominator \(\det(I-yD_s)\) to determine whether it admits an analogous recurrence.
\subsubsection{The recursive relationship between $P_n(y)$ with $Q_n(y)$.}
\begin{defn}
Xin-zhong-zhou\cite{zhou}established a recurrence relation for akins determinants.The $n$-dimensional matrix $J_n$ and $K_n$ be defined as follows for $0\leq i,j\leq n-1$,
\[
J_n=\begin{pmatrix}
2 & -1-y & \ldots & 0 & 0\\
-1 & 2 & \ldots & 0 & 0\\
0 & -1 & \ldots & 0 & 0\\
\vdots & \vdots & \ddots & \vdots & \vdots \\
0 & 0 & \ldots & 2 & -1-y\\
0 & 0 & \ldots & -1 & 2
\end{pmatrix},
K_n=\begin{pmatrix}
2 & -1-y & \ldots & 0 & 0\\
-1 & 2 & \ldots & 0 & 0\\
0 & -1 & \ldots & 0 & 0\\
\vdots & \vdots & \ddots & \vdots & \vdots \\
0 & 0 & \ldots & 2 & -1-y\\
0 & 0 & \ldots & -1 & 1
\end{pmatrix}.
\]

For example,$n=3$ and $n=4$,
\begin{gather*}
J_3=\begin{pmatrix}
2 & -1-y & 0 \\
-1 & 2 & -1-y \\
0 & -1 & 2
\end{pmatrix},
K_3=\begin{pmatrix}
2 & -1-y & 0\\
-1 & 2 & -1-y\\
0 & -1 & 1
\end{pmatrix},\\
J_4=\begin{pmatrix}
2 & -1-y & 0 & 0\\
-1 & 2 & -1-y & 0\\
0 & -1 & 2 & -1-y\\
0 & 0 & -1 & 2
\end{pmatrix},
K_4=\begin{pmatrix}
2 & -1-y & 0 & 0\\
-1 & 2 & -1-y & 0\\
0 & -1 & 2 & -1-y\\
0 & 0 & -1 & 1
\end{pmatrix}.
\end{gather*}
\end{defn}
Let
\[
P_n(y)=det(J_n),\quad Q_n(y)=det(K_n)
\]
as a polynomial with initial value $P_1(y)=2,P_2(y)=3-y$ and $Q_1(y)=1,Q_2(y)=1-y$.Let us assuming $P_0(y)=Q_0(y)=1$.We obtain the following recursive relationships,

\begin{lem}[First-order nonlinear]\label{pqf}
For $n\geq1$,the two polynomials $P_n(y),Q_n(y)$ satisfies the first-order recurrence as follows,
\[
P_{n}(y)=Q_n(y)+P_{n-1}(y),Q_{n}(y)=-yP_{n-1}(y)+(1+y)Q_{n-1}(y).
\]
\end{lem}
\begin{proof}
Observe that the two determinants differ only in their $(n,n)$-entries. We therefore expand the determinant of \(P_n(y)\) along its last row,
\[
(0,0,0,\cdots,-1,2)=(0,0,0,\cdots,-1,1)+(0,0,0,\cdots,0,1),
\]
then
\[
P_n(y)=\underbrace{det\begin{pmatrix}
2 & -1-y & \ldots & 0 & 0\\
-1 & 2 & \ldots & 0 & 0\\
\vdots & \vdots & \ddots & \vdots & \vdots \\
0 & 0 & \ldots & 2 & -1-y\\
0 & 0 & \ldots & -1 & 1
\end{pmatrix}}_{M_1}+\underbrace{det\begin{pmatrix}
2 & -1-y & \ldots & 0 & 0\\
-1 & 2 & \ldots & 0 & 0\\
\vdots & \vdots & \ddots & \vdots & \vdots \\
0 & 0 & \ldots & 2 & -1-y\\
0 & 0 & \ldots & 0 & 1
\end{pmatrix}}_{M_2}
\]
It is no hard to see that $M_1=Q_n(y)$ and $M_2=P_{n-1}(y)$,by expand the last row of $M_2$.Thus
\begin{align}
P_n(y)=Q_n(y)+P_{n-1}(y).\label{P_ny}
\end{align}
While for $Q_n(y)$,we adds the last row of the determinant to second-last row to obtain another determinant,
\[
det\begin{pmatrix}
2 & -1-y & 0 & \ldots & 0 & 0\\
-1 & 2 & -1-y & \ldots & 0 & 0\\
\vdots & \vdots & \vdots & \ddots & \vdots & \vdots \\
0 & 0 & 0 & \ldots & 1 & -y\\
0 & 0 & 0 & \ldots & -1 & 1
\end{pmatrix}
\]
and samely expand the determinant by the last row,
\[
det\begin{pmatrix}
2 & -1-y & 0 & \ldots & 0 \\
-1 & 2 & -1-y & \ldots & 0 \\
\vdots & \vdots & \vdots & \ddots & \vdots \\
0 & 0 & 0 & \ldots & 1
\end{pmatrix}+det\begin{pmatrix}
2 & -1-y & 0 & \ldots & 0 & 0\\
-1 & 2 & -1-y & \ldots & 0 & 0\\
\vdots & \vdots & \vdots & \ddots & \vdots & \vdots \\
0 & 0 & 0 & \ldots & -1 & -y
\end{pmatrix},
\]
consider the last colume of second determinant has only a entry,so expend the determinant by the last colume
\[
-ydet\begin{pmatrix}
2 & -1-y & 0 & \ldots & 0 \\
-1 & 2 & -1-y & \ldots & 0 \\
\vdots & \vdots & \vdots & \ddots & \vdots \\
0 & 0 & 0 & \ldots & 2
\end{pmatrix}_{(n-2)\times(n-2)}.
\]
Therefore having the recursive relation
\begin{align}
Q_n(y)=Q_{n-1}(y)-yP_{n-2}(y),\label{Q_ny}
\end{align}
combining equation \ref{P_ny} with \ref{Q_ny} having
\[
P_{n}(y)=Q_n(y)+P_{n-1}(y),Q_{n+1}(y)=Q_n(y)-yP_{n-1}(y)
\]
then both $P_{n-1}(y)$ are eliminated simultanously
\[
Q_{n+1}(y)+yP_n(y)=(1+y)Q_n(y).
\]
Thus we have the first-order recurrence
\[
P_{n}(y)=Q_n(y)+P_{n-1}(y),Q_{n}(y)=-yP_{n-1}(y)+(1+y)Q_{n-1}(y).
\]
\end{proof}
\begin{lem}[second-order linear]\label{pqs}
For every $n\geq2$,the two polynomials satisfies the second-oeder recursive relationshiops:
\[
P_n(y)=2P_{n-1}(y)-(1+y)P_{n-2}(y),Q_n(y)=2Q_{n-1}(y)-(1+y)Q_{n-2}(y).
\]
\end{lem}
\begin{proof}
We firstly expand by the first colume of the two determinants to obtain two $n-1$-dimensionnal determinaters.
\[
2\underbrace{det\begin{pmatrix}
2 & -1-y & \ldots & 0 & 0\\
-1 & 2 & \ldots & 0 & 0\\
\vdots & \vdots & \ddots & \vdots & \vdots \\
0 & 0 & \ldots & 2 & -1-y\\
0 & 0 & \ldots & -1 & 2
\end{pmatrix}}_{T_1}+
\underbrace{det\begin{pmatrix}
-1-y & 0 & \ldots & 0 & 0\\
-1 & 2 & \ldots & 0 & 0\\
\vdots & \vdots & \ddots & \vdots & \vdots \\
0 & 0 & \ldots & 2 & -1-y\\
0 & 0 & \ldots & -1 & 2
\end{pmatrix}}_{T_2}
\]
Thus $T_1=P_{n-1}(y)$ and for $T_2$,expand the first row of the determinant of first-resulted $n-2$-dimensional determinant.
\[
T_2=-(1+y)det\begin{pmatrix}
2 & -1-y & \ldots & 0 & 0\\
-1 & 2 & \ldots & 0 & 0\\
\vdots & \vdots & \ddots & \vdots & \vdots \\
0 & 0 & \ldots & 2 & -1-y\\
0 & 0 & \ldots & -1 & 2
\end{pmatrix}_{(n-2)\times(n-2)}=-(1+y)P_{n-2}(y).
\]
Thus we obtain the second-order recurrence of $P_n(y),Q_n(y)$,
\[
P_n(y)=2P_{n-1}(y)-(1+y)P_{n-2}(y).
\] 
The recurrence for $Q_{n}(y)$ is obtained in the same way.This complete the lemma.
\end{proof}

Given the initial values and recurrence relations for $P_n(y)$ and \(Q_n(y)\), we obtain the following closed-form expressions for these polynomials.
\begin{lem}\label{pqy}
For $n\geq1$,
\[
P_n(y)=\frac{\alpha^{n+1}-\beta^{n+1}}{\alpha-\beta},Q_n(y)=\frac{\alpha^n+\beta^n}{\alpha+\beta}
\]
where $\alpha=1+i\sqrt{y},\beta=1-i\sqrt{y}$.
\end{lem}
\begin{proof}
We will prove this lemma by induction on $n$.When $n=1$,\[P_1(y)=\frac{\alpha^2-\beta^2}{\alpha-\beta}=\alpha+\beta=2,Q_1(y)=\frac{\alpha+\beta}{\alpha+\beta}=1,\]
satisfying hypothesis clearly.Assuming that lemma holds for $n-1$,but for $n$,by the lemma \ref{pqs},
\begin{align*}
P_n(y)&=2P_{n-1}(y)-(1+y)P_{n-2}(y)=2\frac{\alpha^{n}-\beta^{n}}{\alpha-\beta}-(1+y)\frac{\alpha^{n-1}-\beta^{n-1}}{\alpha-\beta}.
\end{align*}
which derive from subtituting 2 into $\alpha+\beta$ and $1+y$ into $\alpha \beta$.And $Q_n(y)$ can be proved by the smae way.Thus complete prove this closed forms.
\end{proof}
\subsubsection{Calculating the closed-expression of $\Phi_s(y),\Psi_s(y)$.}

\begin{lem}\label{detMs}
For every integer $s\geq0$,
\[
det(I-yD_s)=Q_{s+2}(y).
\]
\end{lem}
\begin{proof}
For convinence setting $M_s=I-yD_s$.We first derive the explicit form of $M_s$ from the definition of $D_s=\bar{B}_sB_s$.Then the form of $M_s$ be showed as index $0\leq i,j\leq s+1$
\[
M_s=\begin{pmatrix}
1-(s+1)y & -sy & -(s-1)y & \ldots & -y & 0 \\
-(s+1)y & 1-sy & -(s-1)y & \ldots & -y & 0 \\
-sy & -sy & 1-(s-1)y & \ldots & -y & 0 \\
\vdots & \vdots & \vdots & \ddots & \vdots & \vdots \\
-2y & -2y & -2y & \vdots & 1-y & 0\\
-y & -y & -y & \vdots & -y & 1
\end{pmatrix}.
\]
At the same time, we apply elementary row and column operations to $M_s$, replacing the first row by the first row minus the second row and the first column by the first column minus the second column.
\[
\begin{pmatrix}
2 & -1 & 0 & \ldots & 0 & 0 \\
-1-y & 1-sy & -(s-1)y & \ldots & -y & 0 \\
0 & -sy & 1-(s-1)y & \ldots & -y & 0 \\
\vdots & \vdots & \vdots & \ddots & \vdots & \vdots \\
0 & -2y & -2y & \vdots & 1-y & 0\\
0 & -y & -y & \vdots & -y & 1
\end{pmatrix}.
\]
And then expand obove matrix determinant by the first row,
\[
2\underbrace{\begin{pmatrix}
1-sy & -(s-1)y & \ldots & -y & 0 \\
-sy & 1-(s-1)y & \ldots & -y & 0 \\
\vdots & \vdots & \ddots & \vdots & \vdots \\
-2y & -2y & \vdots & 1-y & 0\\
-y & -y & \vdots & -y & 1
\end{pmatrix}}_{T_1}+\underbrace{\begin{pmatrix}
-1-y & -(s-1)y & \ldots & -y & 0 \\
0 & 1-(s-1)y & \ldots & -y & 0 \\
\vdots & \vdots & \ddots & \vdots & \vdots \\
0 & -2y & \vdots & 1-y & 0\\
0 & -y & \vdots & -y & 1
\end{pmatrix}}_{T_2}.
\]
Clearly $det(T_1)=det(M_{s-1})$ and $det(T_2)=-(1+y)det(M_{s-2})$,expand by the first colume of $T_2$.Then $det(M_s)$ amazedly satisfies the recurrence
\[
det(M_s)=2det(M_{s-1})-(1+y)det(M_{s-2}).
\]
While for the initial value of $det(M_s)$ by definition,
\begin{gather*}
det(M_0)=det(I-yD_0)=det\begin{pmatrix}
1-y & 0\\
-y & 1
\end{pmatrix}=1-y,\\
det(M_1)=det(I-yD_1)=det\begin{pmatrix}
1-2y & 1-y & 0\\
-2y & 1-y & 0\\
-y & 0 & 1
\end{pmatrix}=1-3y.
\end{gather*}
We find that $\det(M_0)=Q_2(y)$ and $\det(M_1)=Q_3(y)$.Moreover, \(\det(M_s)\) and \(Q_{s+2}(y)\) satisfy the same recurrence relation in \(s\). Therefore, by the initial conditions and the recurrence relation, we conclude that
\[
\det(M_s)=Q_{s+2}(y)
\]
for all \(s\geq 0\).\end{proof}

The following lemma determines the row sums of the adjugate matrix, which will serve as the key ingredient in the subsequent computations of the two numerators.
\begin{lem}\label{uMse}
For all integer $s\geq0$ and each $1\leq i \leq s+2$,
\[
u_s^T(I-yD_s)^*e_i=(1+y)^{s+2-i}Q_i(y).
\]
where $e_i$ denotes the $i$-th standard basis vector of $\mathbb{R}^{s+2}$.
\end{lem}
\begin{proof}
Samely let $M_s=I-yD_s$ and $M_s^*=(M_{ij})$ for $0\leq i,j\leq s+1$,$M_{ij}$ denotes the coordinate $(i,j)$ algebraic complement,then
\[
e_i^T(I-yD_s)^*u_s=e_i^TM_s^*u_s=M_{i0}+M_{i1}+\cdots+M_{i,s+1}.
\]
It is no hard to see that $\sum_{j=0}^{s+1}M_{ij}$ is equal to the determinant that $det(M_s)$ be substituted the $i$-th row all entries into 1,
\[
det\begin{pmatrix}
1-(s+1)y & -sy & \ldots & -(s+3-i)y&-(s+2-i)y&\ldots&-y&0 \\
-(s+1)y & 1-sy & \ldots & -(s+3-i)y&-(s+2-i)y&\ldots&-y&0 \\
-sy & -sy & \ldots & -(s+3-i)y&-(s+2-i)y&\ldots& -y & 0 \\
\vdots & \vdots & \ddots & \vdots &\ldots&\ldots&\ldots& \vdots \\
-(s+4-i)y&-(s+4-i)y&\ldots&1-(s+3-i)y&-(s+2-i)y&\ldots&-y&0\\
1 & 1 & \ldots & 1 & 1 & \ldots & 1 & 1\\
-(s+2-i)y&-(s+2-i)y&\ldots&-(s+2-i)y&-(s+2-i)y&\ldots&-y&0\\
\vdots & \vdots & \ldots & \vdots &\vdots&\ddots&\vdots& \vdots \\
-2y & -2y & \ldots & -2y& -2y&\ldots& 1-y & 0\\
-y & -y & \ldots & -y& -y&\ldots&-y & 1
\end{pmatrix}.
\]
Now we perform elementary following-row operations on the matrix,and the $j$-th row vector denoted by $R_j$ for $i+1\leq j\leq s+2$
\[
R_j+(s+3-j)R_i\rightarrow R_j.
\]
so we can obtain an $(s+2-i)$-dimensional upper triangular sub-matrix as follows,
\[
det\begin{pmatrix}
1-(s+1)y & \ldots & -(s+2-i)y&-(s+1-i)y&\ldots&-y&0 \\
-(s+1)y & \ldots & -(s+2-i)y&-(s+1-i)y&\ldots&-y&0 \\
-sy &\ldots & -(s+2-i)y&-(s+1-i)y&\ldots& -y & 0 \\
\vdots &\ddots & \vdots &\vdots& &\vdots& \vdots \\
-(s+4-i)y&\ldots&-(s+2-i)y&-(s+1-i)y&\ldots&-y&0\\
1 &  \ldots & 1 & 1 & \ldots & 1 & 1\\
0&\ldots&0&1+y&\ldots&(s+1-i)y&(s+2-i)y\\
\vdots &  & \vdots &\vdots&\ddots&\vdots& \vdots \\
0 & \ldots & 0& 0&\ldots& 1+y & 2y\\
0 & \ldots & 0& 0&\ldots&0 & 1+y
\end{pmatrix}.
\]
Then
\[
e_i^TM_s^*u_s=(1+y)^{s+2-i}det\begin{pmatrix}
1-(s+1)y & -sy & \ldots & -(s+2-i)y\\
-(s+1)y & 1-sy & \ldots & -(s+2-i)y\\
-sy & -sy &\ldots & -(s+2-i)y\\
\vdots & \vdots &\ddots \\
-(s+4-i)y&-(s+4-i)y&\ldots&-(s+2-i)y\\
1 & 1 & \ldots & 1 \\
\end{pmatrix}_{i\times i}.
\]
Let us assuming
\[
U_i=\begin{pmatrix}
1-(s+1)y & -sy & \ldots & -(s+2-i)y\\
-(s+1)y & 1-sy & \ldots & -(s+2-i)y\\
-sy & -sy &\ldots & -(s+2-i)y\\
\vdots & \vdots &\ddots &\vdots\\
-(s+4-i)y&-(s+4-i)y&\ldots&-(s+2-i)y\\
1 & 1 & \ldots & 1 \\
\end{pmatrix}_{i\times i}.
\]
Now we perform elementary following-row operations on the above matrix,$C_j$ denotes by the colume vector of $U_i$,
\[
C_j-C_{j+1}\rightarrow C_j,1\leq j\leq i-1.
\]
so we convert $U_i$ to a more familiar form,
\[
U_i\rightarrow\begin{pmatrix}
1-y & -y & -y & \ldots & -y & -(s+2-i)y\\
-1-y & 1-y & -y & \ldots & -y & -(s+2-i)y\\
0 & -1-y & 1-y &\ldots & -y & -(s+2-i)y \\
\vdots & \vdots &\ddots &\vdots &\vdots\\
0&0&0&\ldots&1-y & -(s+2-i)y\\
0 & 0 & 0 &\ldots &0 & 1
\end{pmatrix}_{i\times i}.
\]
Preserving the last row vector has only 1,expand the determinant by the coordinate $(i,i)$,
\[
det(U_i)=det\begin{pmatrix}
1-y & -y & -y & \ldots & -y &-y \\
-1-y & 1-y & -y & \ldots & -y &-y\\
0 & -1-y & 1-y &\ldots & -y &-y\\
\vdots &\vdots& \vdots &\ddots &\vdots &\vdots\\
0&0&0&\ldots&-1-y &1-y
\end{pmatrix}_{(i-1)\times (i-1)}.
\]
We perform elementary following-row operations on the above determinant angin,
\[
C_j-C_{j+1}\rightarrow C_j,1\leq j\leq i-1
\]
then
\[
det(U_i)=det\begin{pmatrix}
2 & -1 & 0 & \ldots & 0 &0 \\
-1-y & 2 & -1 & \ldots &0 &0\\
0 & -1-y & 2 &\ldots & 0 &0\\
\vdots & \vdots&\vdots &\ddots &\vdots &\vdots\\
0&0&0&\ldots&2 &-1\\
0 & 0& 0 & \ldots & -1-y & 1-y
\end{pmatrix}_{(i-1)\times (i-1)},
\]
we can decompose the last colume of the determinant,
\[
(0,0,0,\cdots,-1,1-y)^T=(0,0,0,\cdots,-1,1)^T+(0,0,0,\cdots,0,-y)^T.
\]
gives two following $i-1$-dimensional determinants,
\[
det\begin{pmatrix}
2 & -1 & \ldots & 0 &0 \\
-1-y & 2 & \ldots &0 &0\\
0 & -1-y & \ldots & 0 &0\\
\vdots & \vdots &\ddots &\vdots &\vdots\\
0&0&\ldots&2 &-1\\
0 & 0& \ldots & -1-y & 1
\end{pmatrix}+
det\begin{pmatrix}
2 & -1 & \ldots & 0 &0 \\
-1-y & 2 & \ldots &0 &0\\
0 & -1-y & \ldots & 0 &0\\
\vdots & \vdots&\ddots &\vdots &\vdots\\
0&0&\ldots&2 &0\\
0 & 0& \ldots & -1-y & -y
\end{pmatrix}
\]

and then expand the second determinant by the last entry $-y$,thus
\[
det(U_i)=Q_{i-1}(y)-yP_{i-2}(y).
\]
By the recurrence  \ref{pqf} and \ref{pqs} of $P_n(y)$ and $Q_n(y)$,then
\begin{gather*}
u_s^T(1-yD_s)^*e_i=(1+y)^{s+2-i}Q_i(y).
\end{gather*}
\end{proof}
\begin{lem}\label{uMsu}
For every $s\geq0$,
\[
u_s^T(I-yD_s)^*u_s=P_{s+1}(y).
\]
\end{lem}
\begin{proof}
Decompose the all-one vecter into some standard basis vecter,and by lemma \ref{uMse},
\[
u_s^T(I-yD_s)^*u_s=\sum_{i=1}^{s+2}u_s^T(I-yD_s)^*e_i=\sum_{i=1}^{s+2}(1+y)^{s+2-i}Q_i(y).
\]
Next by lemma \ref{pqy},
\[
Q_i(y)=\frac{\alpha^i+\beta^i}{\alpha+\beta},\alpha=1+i\sqrt{y},\beta=1-i\sqrt{y}.
\]
Substituting the formula of $Q_i(y)$ and $\alpha \beta=1+y$ into the equation,
\[
\sum_{i=1}^{s+2}(1+y)^{s+2-i}\frac{\alpha^i+\beta^i}{\alpha+\beta}=\sum_{i=1}^{s+2}(\alpha \beta)^{s+2-i}\frac{\alpha^i+\beta^i}{\alpha+\beta}.
\]
Using the identity $\sum_{i=1}^nx^{n-i}=(1-x^n)/(1-x)$,we can convert the sum into the follows,
\[
\sum_{i=1}^{s+2}(\alpha \beta)^{s+2-i}\frac{\alpha^i+\beta^i}{\alpha+\beta}=\frac{1}{2}\left(\alpha^{s+2}\frac{1-\beta^{s+2}}{1-\beta}+\beta^{s+2}\frac{1-\alpha^{s+2}}{1-\alpha}\right)=\frac{\alpha^{s+2}-\beta^{s+2}}{\alpha-\beta}.
\]
where $\frac{1}{2}$ be used to replace $\alpha+\beta$ and the last equality derive from $\alpha-\beta=2i\sqrt{y}$.Then this can prove completely
\[
u_s^T(I-yD_s)^*u_s=P_{s+1}(y).
\]
\end{proof}

For another numeretor has the following lemma.
\begin{lem}\label{uMsBu}
For every integer $s\geq1$
\[
u_s^T(I-yD_s)^*\bar{B}_su_s=\frac{(1+y)^{s+2}-Q_{s+2}(y)}{y}.
\]
\end{lem}
\begin{proof}
Perserving the formula,
\[
\bar{B}_su_s=(s+2,s+1,\cdots,1)=\sum_{i=1}^{s+2}(s+3-i)e_i.
\]
Then by lemma \ref{uMse},
\[
u_s^T(I-yD_s)^*\bar{B}_su_s=\sum_{i=1}^{s+2}(s+3-i)u_s^T(I-yD_s)^*e_i=\sum_{i=1}^{s+2}(s+3-i)(1+y)^{s+2-i}Q_i(y).
\]
We konwn the closed-expression of $Q_i(y)=(\alpha^i+\beta^i)/(\alpha+\beta)$ and $\alpha \beta=1+y$,let we assuming $N=s+2$,thus the equality can be converted as follows,
\[
\sum_{i=1}^N(N+1-i)(\alpha \beta)^{N-i}\frac{\alpha^i+\beta^i}{\alpha+\beta}=\frac{1}{2}\sum_{i=1}^N(N+1-i)(\alpha^N\beta^{N-i}+\beta^N\alpha^{N-i})=\frac{\alpha^N+\beta^N-2(\alpha \beta)^N}{-2y}
\]
where the last step uses $(1-\alpha)^2=(1-\beta)^2=-y$ and the identity
\[
\sum_{i=0}^{k-1}(k+1)x^k=\frac{1-(k+1)x^k+kx^{k+1}}{(1-x)^2}.
\]
Then
\[
\frac{\alpha^N+\beta^N-2(\alpha \beta)^N}{-2y}=\frac{(1+y)^{s+2}-Q_{s+2}(y)}{y}.
\]
This complete prove the lemma.
\end{proof}
\begin{proof}[Proof of theorem \ref{pqF}]
Substituting Lemmas \ref{uMsu} and \ref{uMsBu} into Proposition \ref{phipsi}, we obtain the following expressions for $\Phi_s(y)$ and $\Psi_s(y)$.
\[
\Phi_s(y)=\frac{P_{s+1}(y)}{Q_{s+2}(y)},\Psi_s(y)=\frac{(1+y)^{s+2}-Q_{s+2}(y)}{yQ_{s+2}(y)}.
\]
To derive the generating function for $N_d(s)$, we decompose $F(x)$ into the contributions from even and odd dimensions. By definition, the odd-dimensional contribution starts at \(d=3\), and hence can be written as a sum over $d\geq 3$ with \(d\) odd.
\[
\sum_{m\geq1}N_{2m+1}(s)x^{2m+1}=x(\Phi_s(x^2)-N_1(s)).
\]
since $\sum_{m\geq1}N_{2m+1}(s)y^m=\Phi_s(y)-N_1(s)$.And the even-dimensional terms contribute
\[
\sum_{m\geq 0}N_{2m+2}(s)x^{2m+2}=x^2\sum_{m\geq0}N_{2m+2}(s)(x^2)^m=x^2\Psi_s(x^2)
\]
These series consist of the generating function,
\[
F(x)=\sum_{d\geq2}N_d(s)x^d=x(\Phi_s(x^2)-N_1(s))+x^2\Psi_s(x^2).
\]
\end{proof}
\begin{proof}[Proof of theorem \ref{pqR}]
We first derive a recurrence relation for \(\Phi_s(y)\). To this end, using the closed-form expression for \(\Phi_s(y)\) and Lemma \ref{pqf}, we express \(P_{s+1}(y)\) and \(Q_{s+2}(y)\) in the following form,
\[
P_{s+1}(y)=Q_{s+1}(y)-P_s(y),Q_{s+2}(y)=Q_{s+1}(y)-yP_s(y),
\]
then for $s\geq1$,
\[
\Phi_s(y)=\frac{P_{s+1}(y)}{Q_{s+2}(y)}=\frac{Q_{s+1}(y)-P_s(y)}{Q_{s+1}(y)-yP_s(y)}.
\]
Dividing numerator and denominator by $Q_{s+1}(y)$ gives
\[
\Phi_s(y)=\frac{1+\Phi_s(y)}{1-y\Phi_s(y)}.
\]
This proves the first-order recurrence relation for \(\Phi_s(y)\), with initial value $\Phi_0(y)=\frac{2}{1-y}$.We next verify the corresponding first-order recurrence relation for $\Psi_s(y)$. By Theorem \ref{pqF}, the relevant identity can be rewritten in the following equivalent form,
\[
1+y\Psi_s(y)=\frac{(1+y)^{s+2}}{Q_{s+2}(y)}.
\]
then no hard to obtain
\[
\frac{1+y\Psi_s(y)}{1+y\Psi_{s-1}(y)}=\frac{(1+y)Q_{s+1}(y)}{Q_{s+2}(y)}.
\]
Combining $\Phi_s(y)=P_{s+1}(y)/Q_{s+2}(y)$ with lemma \ref{pqf},
\[
\frac{(1+y)Q_{s+1}(y)}{Q_{s+2}(y)}=\frac{Q_{s+2}(y)+yP_{s+1}(y)}{Q_{s+2}(y)}=1+y\Phi_s(y).
\]
This completes the proof that \(\Psi_s(y)\) satisfies the first-order recurrence relation, with initial value $\Psi_0(y)=\frac{3+y}{1-y}$.We now turn to the verification of the second-order recurrence relation. For \(s\geq 2\), applying Lemmas \ref{pqf} and \ref{pqs}, we obtain
\[
P_{s+1}(y)=P_s(y)+Q_s(y)-yP_{s-1}(y),Q_{s+2}(y)=-2P_s(y)+(1+y)Q_s(y).
\]
While for $P_s(y)/Q_{s}(y)$,using lemma \ref{pqf} can convert other form,
\[
\frac{P_s(y)}{Q_s(y)}=\frac{Q_s(y)+P_{s-1}(y)}{Q_s(y)}=1+\Phi_{s-2}(y).
\]
substituting it into $\Phi_s(y)$ and gives the recurrence
\[
\Phi_s(y)=\frac{2+(1-y)\Phi_{s-2}(y)}{1-y-2y\Phi_{s-2}(y)}
\]
with two initial values $\Phi_0(y)=\frac{2}{1-y},\Phi_1(y)=\frac{3-y}{1-3y}$.For $\Psi_s(y)$,by theorem \ref{pqF},
\[
\frac{(1+y)^{s+2}}{1+y\Psi_s(y)}=Q_{s+2}(y).
\]
then by the second-order recurrence of $Q_n(y)$ in lemma \ref{pqs},
\[
\frac{1+y}{1+y\Psi_s(y)}=\frac{2}{1+y\Psi_{s-1}(y)}-\frac{1}{1+y\Psi_{s-2}(y)}
\]
with two initial values $\Psi_0(y)=\frac{3+y}{1-y},\Phi_1(y)=\frac{6+3y+y^2}{1-3y}$.
\end{proof}
\subsubsection{Some concrete example}
We verify the validity of the lattice-point enumeration formula for the corresponding polytopes in several specific cases.
\begin{exam}
For the case of fixed $s=3$,
\begin{gather*}
\Phi_3(y)=\frac{y^2-10y+5}{5y^2-10y+1},
\Psi_3(y)=\frac{y^4+5y^3+10y^2+5y+15}{5y^2-10y+1},\\
F(x)=\frac{x^2(x^8+5x^6-24x^3+5x^2+40x+15)}{5x^4-10x^2+1}.
\end{gather*}
According to the difinition,thus $N_2(3)=15$,$N_4(3)=155$,$N_6(3)=1485$,and $N_3(3)=40,N_5(3)=376$,$N_7(3)=3560$.For the other case of fixed $s=4$,
\begin{gather*}
\Phi_4(y)=\frac{6y^2-20y+6}{-y^3+15y^2-15y+1},
\Psi_4(y)=\frac{y^5+6y^4+15y^3+21y^2+21}{-y^3+15y^2-15y+1},\\
F(x)=-\frac{x^2(x^10+6x^8+15x^6+6x^5+21x^4-84x^3+70x+21)}{x^6-15x^4+15x^2-1}.
\end{gather*}
therefore we get $N_2(4)=21$,$N_4(4)=315$,$N_6(4)=4431$ and $N_3(4)=70$,$N_5(4)=966$,$N_7(4)=13446$.
\end{exam}
\subsection{Circular conditions.}In this subsection, we consider the generating function for the system obtained by adding the constraint $x_d+x_1\leq s+\delta_d$ to the original system. This additional constraint closes the chain of consecutive-coordinate inequalities into a cyclic one, allowing us to compare the resulting generating function with that of the original system. We derive the corresponding functional expressions and recurrence relations. Accordingly, we define
\[
N_s^{cyc}=\sharp\{(x_1,\dots,x_d)\in \mathbb{Z}_{\geq0}^d:x_i+x_{i+1}\leq s+\delta_i,1\leq i\leq d,x_{d+1}=x_1\}.
\]
Be same to theorem \ref{pqF},odd-dimensional series be denoted by $\Phi_s^{cyc}(y)$ and even-dimensional series be denoted by $\Psi_s^{cyc}(y)$,namely
\[
\Phi_s^{cyc}(y)=\sum_{m\geq1}N_{2m+1}^{cyc}(s)y^m,
\Psi_s^{cyc}(y)=\sum_{m\geq1}N_{2m}^{cyc}(s)y^m
\]
and the generating function of $N_d^{cyc}(s)$ be dentoed by $F^{cyc}(x)=\sum_{d\geq2}N_d^{cyc}(s)x^d$.

Then by the transfer matrix method for $m\geq1$,
\begin{gather*}
N_{2m}^{cyc}(s)=\sharp\{(x_1,\dots,x_{2m})\in \mathbb{Z}_{\geq0}^{2m}:x_i+x_{i+1}\leq s+\delta_i,1\leq i\leq 2m\}=tr(D_s^m),\\
N_{2m+1}^{cyc}(s)=\sharp\{(x_1,\dots,x_{2m+1})\in \mathbb{Z}_{\geq0}^{2m+1}:x_i+x_{i+1}\leq s+\delta_i,1\leq i\leq 2m+1\}=tr(D_s^m\bar{B}_s).
\end{gather*}
where $tr(M)$ denotes the \emph{trace} of matrix $M$.Consequently
\begin{gather}
\Phi_s^{cyc}(y)=\sum_{m\geq1}tr(D_s^m\bar{B}_s)y^m=tr((I-yD_s)^{-1}\bar{B}_s)-tr(\bar{B}_s),\label{phiscyc}\\
\Psi_s^{cyc}(y)=\sum_{m\geq1}tr(D_s^m)y^m=tr((I-yD_s)^{-1}-(s+2)\label{psiscyc}.
\end{gather}
By Jacobi's formula, the numerators of the two series admit explicit expressions in $s$. Moreover, by \ref{detMs}, the denominator is given by $\det(I-yD_s)=Q_{s+2}(y)$.Consequently, the two series can be expressed separately in terms of $P_n(y)$ and \(Q_n(y)\).
\begin{lem}\label{trMs}
For every integer $s\geq0$,
\[
\Psi_s^{cyc}(y)=-\frac{y\frac{d}{dy}Q_{s+2}(y)}{Q_{s+2}(y)}.
\]
\end{lem}
\begin{proof}
For any $n\times n$ matrix $A$,let $\{\lambda_1,\dots,\lambda_n\}$ be the set of all eigenvalues of $A$,then $tr(A^m)=\sum_{i=1}^n\lambda_i^m$.Hence
\[
\sum_{m\geq1}tr(A^m)y^m=\sum_{m\geq1}\sum_{i=1}^n\lambda_i^my^m
=\sum_{i=1}^n\frac{\lambda_iy}{1-\lambda_iy}.
\]
Since the charactristic polynomial of $A$ can be written
\[
det(I-yA)=(1-\lambda_1y)\cdots(1-\lambda_ny),
\]
thus we have
\[
\sum_{i=1}^n\frac{\lambda_iy}{1-\lambda_iy}=\frac{-y\frac{d}{dy}det(I-yA)}{det(I-yA)}.
\]
Therefore substituting $D_s$ into this formula and by lemma \ref{detMs},we haved proved the lemma.
\end{proof}
\begin{lem}\label{KN}
For every $n\geq1$ and $0\leq i,j\leq n-1$,
\[
(K_n^{-1})_{ij}=\begin{cases}
\frac{(1+y)^{j-i}P_{i-1}(y)Q_{n-1-j}(y)}{Q_n(y)} &\text{ } i\leq j,\\
\frac{P_{j-1}(y)Q_{n-1-i}(y)}{Q_n(y)} &\text{ } i>j.
\end{cases}
\]
\end{lem}
\begin{proof}
Since $det(K_n)\ne0$ clearly,$K_n$ is reversible.Thus
\[
K_n^{-1}=\frac{K_n^*}{det(K_n)}.
\]
Set $K_{ij}$ as $n-1$-dimensional submatrix of $K_n$,by substracting all entries of $i$-th row and $j$-th colume,then
\[
(K_n^{-1})_{ij}=\frac{(-1)^{i+j}det(K_{ji})}{det(K_n)}.
\]
For $i=j$,$K_{ii}$ can be decomposed into two block matrices,$J_{i-1}$ in the upper-left corner and $K_{n-1-i}$ in the lower-right corner,then $det(K_{ii})=det(J_{i-1})det(K_{n-1-i})$.

And for $0\leq i< j\leq n-1$,
\[
K_{ji}=\begin{pmatrix}
2 & -1-y & 0 & \ldots & 0 & \ldots & 0 & \ldots & 0 & 0\\
-1 & 2 & -1-y & \ldots & 0 & \ldots & 0 & \ldots &0 & 0\\
0 & -1 & 2 & \ldots & 0 & \ldots & 0 & \ldots & 0 & 0\\
\vdots & \vdots & \vdots & \ddots & \vdots &  & \vdots &  & \vdots & \vdots\\
0 & 0 & 0 & \ldots & -1-y & \ldots & 0 & \ldots & 0 & 0\\
\vdots & \vdots & \vdots &  & \vdots & \ddots & \vdots &  & \vdots & \vdots\\
0 & 0 & 0 & \ldots & 0 & \ldots & -1-y & \ldots & 0 & 0\\ 
\vdots & \vdots & \vdots &  & \vdots &  & \ddots &  & \vdots & \vdots\\
0 & 0 & 0 & \ldots & 0 & \ldots & 0 & \ldots & 2 & -1-y\\
0 & 0 & 0 & \ldots & 0 & \ldots & 0 & \ldots & -1 & 1
\end{pmatrix}
\]
Examining the geometric structure of this matrix, we observe that \(J_{i-1}\) appears in the upper-left corner and \(K_{s+1-j}\) in the lower-right corner. The central block is a \((j-i)\)-dimensional lower-triangular square matrix \(O_{j-i}\), whose diagonal entries are all \(-1-y\). Thus, we can perform a block decomposition of the matrix. More precisely, let us write
\[
K_{ji}=\begin{pmatrix}
J_{i-1} & 0 & 0\\
0 & O_{j-i} & 0\\
0 & 0 & K_{n-1-j}
\end{pmatrix}
\]
then
\begin{equation}\label{Kji}
(K_n^{-1})_{ij}=\frac{(-1)^{i+j}det(J_{i-1})det(O_{j-i})det(K_{n-1-j})}{det(K_n)}=\frac{(1+y)^{j-i}det(J_{i-1})det(K_{n-1-j})}{det(K_n)}.
\end{equation}

While for $i>j$,$J_{j-1}$ in the upper-left cornor and $K_{n-1-i}$ in the lower-right cornor,but the square matrix in the center becomes a upper-triangular matrix $\bar{O}_{j-i}$ which all entries are -1 in the diagonal.
\[
K_{ji}=\begin{pmatrix}
J_{j-1} & 0 & 0\\
0 & \bar{O}_{j-i} & 0\\
0 & 0 & K_{n-1-i}
\end{pmatrix}
\]
then
\begin{equation}\label{Kij}
(K_n^{-1})_{ij}=\frac{(-1)^{i+j}det(J_{i-1})det(\bar{O}_{j-i})det(K_{n-1-j})}{det(K_n)}=\frac{det(J_{j-1})det(K_{n-1-i})}{det(K_n)}.
\end{equation}
Therefore by the difinition of $P_n(y)$ and $Q_n(y)$,we can obtain the entry of $K_n^{-1}$.
\end{proof}
\begin{lem}\label{trMsB}
For every integer $s\geq0$,
\[
tr((I-yD_s)^*\bar{B}_s)=\sum_{i=0}^{s+1}U_{s+1-i,i}-\sum_{j=0}^{s}U_{j,s-j}.
\]
where $U_{ij}=(1+y)^{j-i}P_i(y)Q_{s+1-j}(y)$ for $i\leq j$ and $U_{ij}=P_j(y)Q_{s+1-i}(y)$ for $i>j$ and $0\leq i,j\leq s+1$.
\end{lem}
\begin{proof}
For any two square matrices $X$ and $Y$,by jacobi'formula,
\[
\frac{d}{dt}det(X+tY)\Big|_{t=0}=tr(X^*Y).
\]
then respectively substituting $I-yD_s$ and $-\bar{B}_s$ into $X$ and $Y$,which proves the below equation,
\begin{equation}
tr((I-yD_s)^*\bar{B}_s)=-\frac{\partial}{\partial t}det((I-yD_s)-t\bar{B}_s)\Big|_{t=0}.\label{partdet(Ms-B)}
\end{equation}
While $det(\bar{B}_s)\ne0$,then $\bar{B}_s$ is reversible,and for index $0\leq i,j\leq s+1$
\[
(\bar{B}_s^{-1})_{ij}=\begin{cases}
1 &\text{if } i+j=s+1,\\
-1 &\text{if } i+j=s+2,\\
0 &\text{else}.
\end{cases}
\]
thus convert the following formula of \eqref{partdet(Ms-B)} into another form,
\[
det((I-yD_s)-t\bar{B}_s)=det(\bar{B}_s)det(\bar{B}_s^{-1}(I-yD_s)-tI).
\]
Let $W_s(y)=\bar{B}_s^{-1}(I-yD_s)$, subtituting it into equation \eqref{partdet(Ms-B)} and applying jacobi's formula again,
\[
\frac{\partial}{\partial t}det(\bar{B}_s)det(W_s(y)-tI)\Big|_{t=0}=-det(\bar{B}_s)tr(W_s^*(y))=-det(I-yD_s)tr(W_s^{-1}(y)).
\]
Recalls $M_s=I-yD_s$ and
\[
\bar{B}_s^{-1}=\begin{pmatrix}
0 & 0 & 0 & \ldots & 0 & 1\\
0 & 0 & 0 & \ldots & 1 & -1\\
0 & 0 & 0 & \ldots & -1 & 0\\
\vdots & \vdots & \vdots & \ddots & \vdots & \vdots\\
0 & 1 & -1 & \ldots & 0 & 0\\
1 & -1 & 0 & \ldots & 0 & 0
\end{pmatrix},
\]
then
\[
W_s(y)=\begin{pmatrix}
-y & -y & -y & \ldots & -y & 1 \\
-y & -y & -y & \ldots & 1 & -1 \\
-y & -y & -y & \ldots & -1 & 0 \\
\vdots & \vdots & \vdots & \ddots & \vdots & \vdots \\
-y & 1 & -1 & \ldots & 0 & 0\\
1 & -1 & 0 & \ldots & 0 & 0
\end{pmatrix}
\]
Perform elementary row and colume transformations on the matrix as follows,for $0\leq i\leq s$ and $0\leq j\leq [(s+2)/2]$,
\[
R_i-R_{i+1}\rightarrow R_i,C_j\rightarrow C_{s+1-j}.
\]
We encode these elementary row and column operations as matrix transformations. Let \(R_s\) and \(T_s\) denote the corresponding row and column transformation matrices, respectively. Writting this matrices,
\[
R_s=\begin{pmatrix}
1 & -1 & 0 & \ldots & 0 & 0\\
0 & 1 & -1 & \ldots & 0 & 0\\
0 & 0 & 1 & \ldots & 0 & 0\\
\vdots & \vdots & \vdots & \ddots & \vdots & \vdots\\
0 & 0 & 0 & \ldots & 1 & -1\\
0 & 0 & 0 & \ldots & 0 & 1
\end{pmatrix},
T_s=\begin{pmatrix}
0 & 0 & 0 & \ldots & 0 & 1\\
0 & 0 & 0 & \ldots & 1 & 0\\
0 & 0 & 0 & \ldots & 0 & 0\\
\vdots & \vdots & \vdots & \ddots & \vdots & \vdots\\
0 & 1 & 0 & \ldots & 0 & 0\\
1 & 0 & 0 & \ldots & 0 & 0
\end{pmatrix}.
\]
Then the transformed matrix is given by
\[
R_sW_s(y)T_s=\begin{pmatrix}
2 & -1-y & 0 & \ldots & 0 & 0\\
-1 & 2 & -1-y & \ldots & 0 & 0\\
0 & -1 & 2 & \ldots & 0 & 0\\
\vdots & \vdots & \vdots & \ddots & \vdots & \vdots\\
0 & 0 & 0 & \ldots & 2 & -1-y\\
0 & 0 & 0 & \ldots & -1 & 1
\end{pmatrix}=K_{s+2}.
\]
Clearly $R_s$ and $T_s$ are reversible,
\[
W_s^{-1}(y)=T_sK_{s+2}^{-1}R_s.
\]
By lemma \ref{KN} and setting $n=s+2$,then for $0\leq i,j\leq s+1$,
\[
(K_{s+2}^{-1})_{ij}=\begin{cases}
\frac{(1+y)^{j-i}P_i(y)Q_{s+1-j}(y)}{Q_{s+2}(y)} &\text{ } i\leq j,\\
\frac{P_j(y)Q_{s+1-i}(y)}{Q_{s+2}(y)} &\text{ } i>j.
\end{cases}
\]
where by lemma \ref{pqy},assuming $P_0(y)=Q_0(y)=1$.Thus
\[
tr(W_s^{-1}(y))=tr(T_sK_{s+2}^{-1}R_s)=\sum_{i=0}^{s+1}(K_{s+2}^{-1})_{s+1-i,i}-\sum_{j=0}^{s}(K_{s+2}^{-1})_{j,s-j}.
\]
Let $U_{ij}=Q_{s+2}(y)(K_{s+2}^{-1})_{ij}$,by lemma \ref{detMs},
\[
Q_{s+2}(y)tr(W_s^{-1}(y))=\sum_{i=0}^{s+1}U_{s+1-i,i}-\sum_{j=0}^{s}U_{j,s-j}.
\]
Therefore we had proved
\[
tr((I-yD_s)^*\bar{B}_s)=\sum_{i=0}^{s+1}U_{s+1-i,i}-\sum_{j=0}^{s}U_{j,s-j}.
\]
\end{proof}
\begin{proof}[Proof of theorem \ref{pqC}]
Let us examine the realtionship between the derivative of $Q_n(y)$ and $P_n(y)$,by lemma \ref{pqf},
\[
Q_{s+2}(y)=\frac{(1+i\sqrt{y})^{s+2}+(1-i\sqrt{y})^{s+2}}{2},
\]
then deriving the derivative to $y$
\[
\frac{d}{dy}Q_{s+2}(y)=-\frac{s+2}{2}\frac{(1+i\sqrt{y})^{s+1}-(1-i\sqrt{y})^{s+1}}{2i\sqrt{y}}=-\frac{s+2}{2}P_s(y).
\]
Thus by lemma \ref{trMs},we can prove that
\[
\Psi_s^{cyc}(y)=\frac{y(s+2)P_s(y)}{2Q_{s+2}(y)}.
\]
And by equation \eqref{phiscyc},we have the form of adjugate matrix
\[
\Phi_s^{cyc}(y)=\frac{tr((I-yD_s)^*\bar{B}_s)}{det(I-yD_s)}-tr(\bar{B}_s).
\]
Therefore by lemma \ref{detMs} and \ref{trMsB},
\[
\Phi_s^{cyc}(y)=\frac{\sum_{i=0}^{s+1}U_{s+1-i,i}-\sum_{j=0}^{s}U_{j,s-j}}{Q_{s+2}(y)}-\left[\frac{s+3}{2}\right].
\]
Recalls the definition of the generating function $F^{cyc}(x)$,decomposing the sums of odd- and even- terms,
\[
F^{cyc}(x)=\sum_{m\geq1}N_{2m}^{cyc}(s)x^{2m}+\sum_{m\geq1}N_{2m+1}^{cyc}(s)x^{2m+1}=\Psi_s^{cyc}(x^2)+x\Phi_s^{cyc}(x^2).
\]
This complete prove the theorem.
\end{proof}
\subsubsection{Cyclic concrete example}
\begin{exam}
For the case of fixed $s=3$,
\begin{gather*}
\Phi_3^{cyc}(y)=\frac{y^4+3y^3-7y^2+33y}{5y^2-10y+1},\Psi_3^{cyc}(y)=\frac{-10y^2+10y}{5y^2-10y+1},\\
F^{cyc}(x)=\frac{x^9+3x^7-7x^5-10x^4+33x^3+10x^2}{4x^4-10x^2+1},
\end{gather*}
thus $N_2^{cyc}(3)=10$,$N_4^{cyc}(3)=90$,$N_6^{cyc}(3)=850$ and $N_3^{cyc}(3)=33$,$N_5^{cyc}(3)=323$,$N_7^{cyc}(3)=3068$.While for the case of fixed $s=4$,
\begin{gather*}
\Phi_4^{cyc}(y)=\frac{-y^5-4y^4-14y^2+38y^2-57y}{y^3-15y^2+15y-1},\Psi_4^{cyc}(y)=\frac{-3y^3+10y^2-5y}{y^3-15y^2+15y-1},\\
F^{cyc}(x)=\frac{-x^11-4x^9-14x^7-3x^6+38x^5+30x^4-57x^3-15x^2}{(x^4-14x^2+1)(x^2-1)},
\end{gather*}
therefore $N_2^{cyc}(4)=15$,$N_4^{cyc}(4)=195$,$N_6^{cyc}(4)=2703$ and $N_3^{cyc}(4)=57$,$N_5^{cyc}(4)=817$,$N_7^{cyc}(4)=11414$.
\end{exam}
\begin{flushleft}
\textbf{Remark:}\emph{Gives two polynomials $\hat{P}_n(y)$ and $\hat{Q}_n(y)$ as the odd and even parts of $(1+\sqrt{y})^n$,thus
\[
\hat{P}_n(y)=\sum_{k=0}^{[(n-1)/2]}\binom{n}{2k+1}y^k,\hat{Q}_n(y)=\sum_{k=0}^{[n/2]}\binom{n}{2k}y^k.
\]
Since Pascal's identity
\[
\binom{n}{2k+1}=\binom{n-1}{2k+1}+\binom{n-1}{2k},
\]
gives the second-order recursive relationship for $n\geq2$
\[
\hat{P}_n(-y)=2\hat{P}_{n-1}(-y)-(1+y)\hat{P}_{n-2}(-y),\hat{Q}_n(-y)=2\hat{Q}_{n-1}(-y)-(1+y)\hat{Q}_{n-2}(-y).
\]
But the initial value $\hat{P}_1(-y)=1,\hat{P}_2(-y)=2$ and $\hat{Q}_1(-y)=1,\hat{Q}_2(-y)=1-y$,then for $n\geq1$
\[
\hat{P}_n(y)=P_{n-1}(y),\hat{Q}_n(-y)=Q_n(y).
\]
Hence we can also denote the all above conclusions by $\hat{P}_n(-y)$ and $\hat{Q}_n(-y)$.}
\end{flushleft}
We change the variable $y$ as $\theta$,set $y=\tan^2 \theta$,so $\sqrt{y}=\tan \theta$ and $\theta=\arctan y$.Since the identity $1+i\tan \theta=e^{i\theta}/\cos \theta$,
\[
\alpha=\frac{e^{i\theta}}{\cos \theta},\beta=\frac{e^{-i\theta}}{\cos \theta}.
\]
Thus by lemma \ref{pqy} gives the following trigonometric closed forms
\begin{equation}
P_n(y)=\frac{\sin(n+1)\theta}{\sin \theta \cos^n\theta},Q_n(y)=\frac{\cos n\theta}{\cos^n\theta}.\label{tripq}
\end{equation}
\begin{lem}[trigonometric forms]\label{triPQ}
For every integer $s\geq1$,
\begin{gather*}
\Phi_s(y)=\frac{\tan((s+2)\arctan \sqrt{y})}{\sqrt{y}},\\
\Psi_s(y)=\frac{1}{y}\left(\frac{1}{\cos((s+2)\arctan \sqrt{y})\cos^{s+2}(\arctan \sqrt{y})}-1\right).
\end{gather*}
\end{lem}
\begin{proof}
By theorem \ref{pqF} and equation \eqref{tripq},
\[
\Phi_s(y)=\frac{P_{s+1}(y)}{Q_{s+2}(y)}=\frac{\tan(s+2)\theta}{\tan \theta}=\frac{\tan((s+2)\arctan \sqrt{y})}{\sqrt{y}}.
\]
Samely for $\Psi_s(y)$,gives
\[
\Psi_s(y)=\frac{(1+y)^{s+2}-Q_{s+2}(y)}{yQ_{s+2}(y)}=\frac{1}{y}\left(\frac{1}{\cos(s+2)\theta \cos^{s+2}\theta}-1\right),
\]
substituting $\theta=\arctan \sqrt{y}$ yields
\[
\Psi_s(y)=\frac{1}{y}\left(\frac{1}{\cos((s+2)\arctan \sqrt{y})\cos^{s+2}(\arctan \sqrt{y})}-1\right).
\]
\end{proof}
\section{Ehrhart series and $h^*$-vectors}
In this section, we return to the fixed-dimensional Ehrhart theory of the polytope \(P_d^{(s)}\). Unlike in Section 2, where \(s\) is fixed and \(d\) varies, we now fix the dimension \(d\) and vary \(s\). In this setting, we investigate the lattice-point enumeration of the dilates of the polytope and study the resulting Ehrhart polynomial and Ehrhart series.
\subsection{Ehrhart framwork and some amazed formula }.Recall that the lattice-point enumerator of the dilates of \(\mathcal{P}_d^{(s)}\) is given by the Ehrhart polynomial. By Ehrhart theory, for every fixed integer \(d\geq 2\) and \(s\geq 1\), this lattice-point enumerator is a polynomial in the dilation parameter,
\[
Ehr_{\mathcal{P}_s^{(s)}}(z)=1+\sum_{m\geq1}L_d^{(s)}(m)z^m=\frac{h_d^{(s)}(z)}{(1-z)^{d+1}},L_d^{(s)}(m)=\sharp(m\mathcal{P}_s^{(s)}\cap \mathbb{Z}^d).
\]
defining the $h^*$-polynomial $h_d^{(s)}(z)=\sum_{k=0}^dh_{d,k}^{(s)}z^k$.Thus $h_{d,k}^{(s)}$ denotes the $k$-th Ehrhart $h^*$-coefficient.Difine the following generating function
\[
G_m^{(s)}:=\sum_{d\geq0}L_d^{(s)}(m)y^d,H_k^{(s)}(y):=\sum_{d\geq0}h_{d,k}^{(s)}y^d.
\]
where for $d=0$ and $d=1$,set
\[
\mathcal{P}_0^{(s)}=\{0\},\mathcal{P}_1^{(s)}=\{x\in \mathbb{R}:0\leq x\leq s+1\},
\]
then gives
\[
L_0^{(s)}(m)=1,L_1^{(s)}(m)=(s+1)m+1.
\]

\begin{lem}\label{Hksy}
For fixed $s\geq1$ and $k\geq0$,the generating function of the $k$-th Ehrhart $h^*$-coefficient as follows
\[
H_k^{(s)}(y)=\sum_{i=0}^k\frac{(-1)^i}{i!}y^{i-1}\frac{d^i}{dy^i}(yG_{k-i}^{(s)}(y)).
\]
\end{lem}
\begin{proof}
Recalls that
\[
1+\sum_{m\geq1}L_d^{(s)}(m)z^m=\sum_{m\geq0}L_d^{(s)}(m)z^m=\frac{h_d^{(s)}z}{(1-z)^{d+1}}=\frac{\sum_{k=0}^dh_{d,k}^{(s)}z^k}{(1-z)^{d+1}},
\]
then
\[
\sum_{k=0}^dh_{d,k}^{(s)}z^k=(1-z)^{d+1}\left(\sum_{m\geq0}L_d^{(s)}(m)z^m\right).
\]
Since $(1-z)^{d+1}=\sum_{i=0}^{d+1}\binom{d+1}{i}(-z)^i$,
\[
h_{d,k}^{(s)}=\sum_{i=0}^{k}(-1)^i\binom{d+1}{i}L_d^{(s)}(m).
\]
For each fixed $i\geq0$,we convert the definition of $G_m^{(s)}(y)$,
\[
yG_{k-i}^{(s)}(y)=\sum_{d\geq0}L_d^{(s)}(k-i)y^{d+1}
\]
then the $i$-th derivative of the expression with respect to $y$,
\[
\frac{d^i}{dy^i}(yG_{k-i}^{(s)}(y))=\sum_{d\geq0}i!\binom{d+1}{i}L_d^{(s)}(k-i)y^{d+1-i},
\]
Multiply both sides by $y^{i-1}/i!$ simultaneously and by the definition of $h^*$-polynomial to get the lemma.
\end{proof}
\subsection{the dimension generating function and explicit exprssion}Lemma \ref{Hksy} implies that $H_k^{(s)}(y)$ can be expressed in terms of the series \(G_m^{(s)}(y)\). We therefore proceed to determine explicit expressions for \(G_m^{(s)}(y)\) for all integers \(s,m\geq 1\). To this end, we define two transfer matrices, $X_m^{(s)}$ and $Y_m^{(s)}$, indexed by $0\leq i,j\leq (s+1)m$,
\[
(X_m^{(s)})_{i,j}=\begin{cases}
1,&\text{if } i+j\leq (s+1)m,\\
0,&\text{else }.
\end{cases},
(Y_m^{(s)})_{i,j}=\begin{cases}
1,&\text{if } i+j\leq sm,\\
0,&\text{else }.
\end{cases}.
\]
Set $Z_m^{(s)}=X_m^{(s)}Y_m^{(s)}$ and all-one $(s+1)m+1$-dimensional colume vector be denoted by $u$.Decompose $G_m^{(s)}(y)$ into odd and even parts
\[
G_m^{(s)}(y)=\sum_{d\geq0}L_d^{(s)}(m)y^d=\sum_{r\geq0}L_{2r}^{(s)}(m)y^{2r}+\sum_{r\geq0}L_{2r+1}^{(s)}y^{2r+1}.
\]
Let $U_m^{(s)}(z)$ as the even parts and $V_m^{(s)}(z)$ as the odd part,
\[
U_m^{(s)}(z)=\sum_{r\geq1}L_{2r}^{(s)}(m)z^{r-1},V_m^{(s)}(z)=\sum_{r\geq0}L_{2r+1}^{(s)}(m)z^{r}.
\]
then by transfer matrix theory,
\begin{gather}
U_m^{(s)}(z)=\sum_{r\geq1}u^T(Z_m^{(s)})^{r-1}X_m^{(s)}uz^{r-1}=u^T(I-zZ_m^{(s)})^{-1}X_m^{(s)}u,\label{Umsz}\\
V_m^{(s)}(z)=\sum_{r\geq0}u^T(Z_m^{(s)})^ruz^r=u^T(I-zZ_m^{(s)})^{-1}u.\label{Vmsz}
\end{gather}
where $I$ is $((s+1)m+1)$-dimensional identity matrix.

\begin{lem}\label{UVmsz}
Let $J_m^{(s)}$ denotes $(sm+1)$-dimensional square matrix,
\[
(J_m^{(s)})_{ij}=\begin{cases}
1 &\text{if \quad} i+j\leq sm,\\
0 &\text{else}.\end{cases}
\]
and let $R_m^{(s)}$ denotes block matrix $(J_m^{(s)}\quad 0)\in \mathbb{R}^{(sm+1)\times ((s+1)m+1)}$,and another block matrix $I_m^{(s)}=(I_{sm+1}\quad 0)^T\in \mathbb{R}^{((s+1)m+1)\times(sm+1)}$.Then
\begin{gather}
U_m^{(s)}(z)=\frac{((s+1)m+2)((s+1)m+1)}{2}+z\mathbf{\mu}^T(I-zR_m^{(s)}T_m^{(s)})^{-1}\mathbf{\omega},\label{umsz}\\
V_m^{(s)}(z)=(s+1)m+1+z\mathbf{\mu}^T(I-zR_m^{(s)}T_m^{(s)})^{-1}\mathbf{\nu}.\label{vmsz}
\end{gather}
where $\mathbf{\mu}^T=((s+1)m+1,(s+1)m,\cdots,m+1)$,$\mathbf{\nu}=(sm+1,sm,\cdots,1)^T$ and $\mathbf{\omega}=R_m^{(s)}X_m^{(s)}u$.
\end{lem}
\begin{proof}
Preserving the transfer matrix $Y_m^{(s)}$,we find
\[
Y_m^{(s)}=\begin{pmatrix}
J_m^{(s)} & 0\\
0 & 0
\end{pmatrix}=I_m^{(s)}R_m^{(s)}.
\]
Thus $Z_m^{(s)}=X_m^{(s)}Y_m^{(s)}=X_m^{(s)}I_m^{(s)}R_m^{(s)}$.Let $T_m^{(s)}=X_m^{(s)}I_m^{(s)}$ and then $Z_m^{(s)}=T_m^{(s)}R_m^{(s)}$.

We record the following push-through identity, which will be used in the subsequent compression argument. Let $A$ and \(B\) be rectangular matrices of compatible dimensions, and assume that $I-zBA$ is invertible. Setting $C=(I-zBA)^{-1}$,we obtain
\[
(I-zBA)C=C-zBAC=I,
\]
then 
\[
(I-zAB)(I+zACB)=I+zACB-zAB-z^2ABACB=I+zA(C-I-zBAC)B=I
\]
and one has the formal poweer-series identity
\begin{equation}
(I-zAB)^{-1}=I+zA(I-zBA)^{-1}B.\label{IABN}
\end{equation}
and multiply on both sides $A$ by the right
\begin{equation}
(I-zAB)^{-1}A=A(I+z(I-zBA)^{-1}BA)=A(I-zBA)^{-1}.\label{IABNA}
\end{equation}
Substituting $A=T_m^{(s)}$ and $B=R_m^{(s)}$ into \ref{IABN},gives
\begin{equation}
(I-zT_m^{(s)}R_m^{(s)})^{-1}=(I-zZ_m^{(s)})^{-1}=I+zT_m^{(s)}(I-zR_m^{(s)}T_m^{(s)})^{-1}R_m^{(s)},\label{ITRN}
\end{equation}
where the $(i,j)$-entry of $R_m^{(s)}T_m^{(s)}$ for $0\leq i,j\leq sm$,
\[
(R_m^{(s)}T_m^{(s)})_{ij}=\sum_{k=0}^{(s+1)m}(R_m^{(s)})_{ik}(T_m^{(s)})_{kj},
\]
but the entry of $R_m^{(s)}$ and $T_m^{(s)}$ are either 0 or 1,thus
\[
\sum_{k=0}^{(s+1)m}(R_m^{(s)})_{ik}(T_m^{(s)})_{kj}=\sharp\{i+k\leq sm,k+j\leq (s+1)m,k\geq0\},
\]
then
\[
(R_m^{(s)}T_m^{(s)})_{ij}=\min\{sm-i,(s+1)m-j\}+1.\label{RTmsij}
\]

Recalls $V_m^{(s)}(z)=u^T(I-zZ_m^{(s)})^{-1}u$,and by equation \ref{ITRN},
\[
V_m^{(s)}(z)=(s+1)m+1+zu^TT_m^{(s)}(I-zR_m^{(s)}T_m^{(s)})^{-1}R_m^{(s)}u
\]
Since $T_m^{(s)}=X_m^{(s)}I_m^{(s)}$,let
\begin{gather*}
\mathbf{\mu}:=u^TT_m^{(s)}=((s+1)m+1,(s+1)m,\cdots,m+1)^T\in \mathbb{R}^{sm+1},\\
\mathbf{\nu}:=R_m^{(s)}u=(sm+1,sm,\cdots,1)^T\in \mathbb{R}^{sm+1}.
\end{gather*}
then
\[
V_m^{(s)}(z)=(s+1)m+1+z\mathbf{\mu}^T(I-zR_m^{(s)}T_m^{(s)})^{-1}\mathbf{\nu}
\]
Combining equation \ref{IABN} with \ref{IABNA},since $Y_m^{(s)}=I_m^{(s)}R_m^{(s)}$ and $T_m^{(s)}=X_m^{(s)}I_m^{(s)}$,we get
\[
U_m^{(s)}(z)=u^TX_m^{(s)}(I-zI_m^{(s)}R_m^{(s)}X_m^{(s)})^{-1}u=u^TX_m^{(s)}u+zu^TT_m^{(s)}(I-zR_m^{(s)}T_m^{(s)})^{-1}R_m^{(s)}X_m^{(s)}u
\]
Since $\mathbf{\mu}=u^TT_m^{(s)}$ and $u^TX_m^{(s)}u=(sm+m+2)(sm+m+1)/2$,let $\mathbf{\omega}=R_m^{(s)}X_m^{(s)}u\in \mathbb{R}^{sm+1}$,and for $0\leq i\leq sm$,$0\leq j\leq (s+1)m$,
\[
(R_m^{(s)}X_m^{(s)})_{ij}=\sharp \{i+k\leq sm,k+j\leq (s+1)m\}=sm-\max\{i,j-m\},
\]
then
\[
U_m^{(s)}(z)=\frac{((s+1)m+2)((s+1)m+1)}{2}+z\mathbf{\mu}^T(I-zR_m^{(s)}T_m^{(s)})^{-1}\mathbf{\omega}.
\]
which complete vertify \ref{umsz} and \ref{vmsz}.
\end{proof}

As in Section 2, we decompose the generating function into its odd- and even-dimensional parts. To obtain explicit expressions for the two components, we introduce several auxiliary quantities and identities, which will be used to simplify and evaluate the resulting formulas,for  $0\leq j\leq s$
\begin{gather*}
\bar{\alpha}_j^{(s,m)}:=\binom{(s+1-j)m+j+1}{2j},\quad \bar{\beta}_j^{(s,m)}=\binom{(s+1-j)m+j+1}{2j+1},\\
\bar{\rho}_j^{(s,m)}=\bar{\alpha}_j^{(s,m)}
+\bar{\beta}_j^{(s,m)}=\binom{(s-j+1)m+j+2}{2j+1}.
\end{gather*}

Define the $k$-th component of two families vector $\sigma_j,\tau_j\in \mathbb{R}^{sm+1}$ by
\[
\sigma_j(k)=\binom{(k-jm)_++j}{2j},\quad \tau_j(k)=\binom{(k-jm)_++j+1}{2j+1}.
\]
where the symbol $(x)_+=\max \{x,0\}$.

Setting $C_m^{(s)}=R_m^{(s)}T_m^{(s)}\in \mathbb{R}^{(sm+1)\times(sm+1)}$ and recalls $u_s\in \mathbb{R}^{s+2}$ is all-one colume vector.
\begin{lem}\label{sigmatau}
For $0\leq j\leq s-1$,gives the identity
\begin{gather}
\sigma_j^TC_m^{(s)}=\bar{\alpha}_{j+1}^{(s,m)}u_{sm-1}^T-\sigma_{j+1}^T,\label{sigmaC}\\
\tau_j^TC_m^{(s)}=\bar{\rho}_{j+1}^{(s,m)}u_{sm-1}^T-\tau_{j+1}^T.\label{tauC}
\end{gather}
where $\bar{\rho}_{j+1}^{(s,m)}=\bar{\alpha}_{j+1}^{(s,m)}+
\bar{\beta}_{j+1}^{(s,m)}$.
\end{lem}
\begin{proof}
Since $k\leq sm$,then for $\max \{k-sm,0\}=0$ and $\sigma_s=\tau_s=\mathbf{0}$.By lemma \ref{UVmsz} for $0\leq a,b\leq sm$,
\[
(C_m^{(s)})_{ab}=\min \{sm-a,sm+m-b\}+1=sm-a+1-(b-a-m)_+,
\]
then for fixed $j=0$,thus $\sigma_0=\mathbf{1}$,
\[
(\sigma_0^TC_m^{(s)})(k)=\sum_{i=0}^{sm}(sm-i+1-(k-i-m)_+)=\begin{cases}
\frac{(sm+2)(sm+1)}{2} &\text{}\quad k\leq m,\\
\frac{((s+1)m-k+1)(sm+k+2)}{2}&\text{}\quad k>m.
\end{cases}
\]
and $\bar{\alpha}_1^{(s,m)}-\sigma_1(k)=\frac{((sm+2)(sm+1)-((k-m)_++1)((k-m)_+))}{2}$,then $\sigma_0^TC_m^{(s)}=\bar{\alpha}_1u_{sm+1}^T-\sigma_1^T$.While for $j\geq 1$,$\binom{j}{2j}=0$,then
\begin{align}
(\sigma_j^TC_m^{(s)})(k)&=\sum_{i=0}^{sm}\sigma_j(i)(C_m^{(s)})_{ik}=\sum_{i=0}^{sm}(sm-i+1-(k-i-m)_+)\binom{(i-jm)_++j}{2j}\notag\\
&=\sum_{i=jm}^{sm}(sm-i+1-(k-i-m)_+)\binom{i-jm+j}{2j}.\label{sigmaCmsk1}
\end{align}
Let us assuming $i-jm=t$,and  by the hockey-stick identity 
\[
\sum_{t=0}^{N}(N-t+1)\binom{t+j}{2j}=\binom{N+j+2}{2j+2},
\]
substitutes $N=(s-j)m$ into the above formula,
\[
\binom{(s-j)m+j+2}{2j+2}=\bar{\alpha}_{j+1}^{(s,m)}.
\]
then by equation \ref{sigmaCmsk1},
\begin{equation}
(\sigma_j^TC_m^{(s)})(k)=\bar{\alpha}_{j+1}^{(s,m)}-\sum_{t=0}^{(s-j)m}(k-(j+1)m-t)_+\binom{t+j}{2j}.\label{sigmaCmsk2}
\end{equation}
In the rear sum part,when $k\leq (j+1)m$,$(k-(j+1)m-t)_+=0$ and $\sigma_{j+1}(k)=0$ for every $k\leq sm$.Then for $k> (j+1)m$,
\begin{align*}
\sum_{t=0}^{(s-j)m}(k-(j+1)m-t)_+\binom{t+j}{2j}=\sum_{t=0}^{k-(j+1)m-1}(k-(j+1)m-t)\binom{t+j}{2j}.
\end{align*}
By the hockey-stick identity again,
\[
\sum_{t=0}^{k-(j+1)m-1}(k-(j+1)m-t)\binom{t+j}{2j}=\binom{k-(j+1)m+1+j}{2j+2}=\sigma_{j+1}(k).\label{sigmaCmsk3}
\]
Thus by \ref{sigmaCmsk2} and \ref{sigmaCmsk3},for $0\leq k\leq sm$
\[
(\sigma_j^TC_m^{(s)})(k)=\bar{\alpha}_{j+1}^{(s,m)}-\sigma_{j+1}^T(k)
\]
then for $0\leq j\leq s-1$,we get the identity $\sigma_j^TC_m^{(s)}=\bar{\alpha}_{j+1}^{(s,m)}u_{sm-1}^T-\sigma_{j+1}^T$.

We can prove the second identity \ref{tauC} by the same way,we fixed $j=0$ first,then $\tau_0(k)=k+1$,
\[
(\tau_0^TC_m^{(s)})(k)=\sum_{i=0}^{sm}(i+1)(sm-i+1-(k-i-m)_+).
\]
If $k\leq m$,$(k-i-m)_+=0$ and $(\tau_0^TC_m^{(s)})(k)=\binom{sm+3}{3}$.If $k>m$,
\[
(\tau_0^TC_m^{(s)})(k)=\sum_{i=0}^{k-m}(i+1)((s+1)m-k+1)=\frac{((s+1)m-k+1)(k-m+2)(k-m+1)}{2}
\]
While
\[
(\bar{\alpha}_1^{(s,m)}+\bar{\beta}_1^{(s,m)})-\tau_1(k)=\begin{cases}
\binom{sm+3}{3}&\text{}\quad k\geq m,\\
\binom{sm+3}{3}-\binom{k-m+2}{3}&\text{}\quad k>m.
\end{cases}
\]
then $\tau_0^TC_m^{(s)}=(\bar{\alpha}_1^{(s,m)}+\bar{\beta}_1^{(s,m)})u_{sm-1}^T-\tau_{1}^T$.
For $j\geq1$,$\binom{j+1}{2j+1}=0$ and
\begin{align}
(\tau_j^TC_m^{(s)})(k)&=\sum_{i=0}^{sm}(sm-i+1-(k-i-m)_+)\binom{(i-jm)_++j+1}{2j+1}\notag \\
&=\sum_{i=jm}^{sm}(sm-i+1-(k-i-m)_+)\binom{i-jm+j+1}{2j+1}.\label{tauCmsk1}
\end{align}
Let $i-jm=t$,convert \ref{tauCmsk1} into two parts form,
\[
\sum_{t=0}^{(s-j)m}((s-j)m-t+1)\binom{t+j+1}{2j+1}-\sum_{t=0}^{(s-j)m}((k-(j+1)m-t)_+)\binom{t+j+1}{2j+1}.
\]
While for the first part,let $N=(s-j)m$
\[
\sum_{t=0}^{N}(N-t+1)\binom{t+j+1}{2j+1}=\binom{N+j+3}{2j+3},
\]
by the Pascal's identity,
\begin{equation}
\binom{N+j+3}{2j+3}=\binom{N+j+2}{2j+2}+\binom{N+j+2}{j+3}=\bar{\alpha}_{j+1}^{(s,m)}+\bar{\beta}_{j+1}^{(s,m)}.\label{tauCmsk2}
\end{equation}
In the rear sum part,when $k\leq (j+1)m$,$(k-(j+1)m-t)_+=0$ and $\tau_{j+1}(k)=0$ for every $k\leq sm$.Then for $k> (j+1)m$,
\begin{align}
\sum_{t=0}^{(s-j)m}((k-(j+1)m-t)_+)\binom{t+j+1}{2j+1}&=\sum_{t=0}^{k-(j+1)m-1}(k-(j+1)m-t)\binom{t+j+1}{2j+1}\notag \\
&=\binom{k-(j+1)m+j+2}{2j+3}=\tau_{j+1}(k).\label{tauCmsk3}
\end{align}
Therefore by \ref{tauCmsk2} and \ref{tauCmsk3},for $0\leq j\leq s-1$
\[
\tau_j^TC_m^{(s)}=\bar{\rho}_{j+1}^{(s,m)}
u_{sm+1}^T-\tau_{j+1}^T.
\]
which have proved this identity of lemma \ref{sigmatau}
\end{proof}
\begin{lem}\label{omega}
For $0\leq j\leq s-1$,
\begin{gather*}
\sigma_j^T\mathbf{\omega}=\eta_{j+1}^{(s,m)}\bar{\alpha}_{j+1}^{(s,m)},\tau_j^T\mathbf{\omega}=\xi_{j+1}^{(s,m)}\bar{\rho}_{j+1}^{(s,m)},\\
\sigma_j^T\mathbf{\nu}=\bar{\alpha}_{j+1}^{(s,m)},\tau_j^T\mathbf{\nu}=\bar{\rho}_{j+1}^{(s,m)}.
\end{gather*}
where $\eta_j^{(s,m)}=\frac{j((2s+3)m+3)}{2j+1}$ and $\xi_j^{(s,m)}=\frac{((2j+1)s+3j+1)m+3j+1}{2j+2}$.
\end{lem}
\begin{proof}
Recalls $\mathbf{\omega}=R_m^{(s)}X_m^{(s)}u$,where $u$ is a $((s+1)m+1)$-dimensional all-one colume vector.We konw $\sigma_s=\tau_s=0$,then $\sigma_s^T\mathbf{\omega}=\tau_s^T\mathbf{\omega}=0$,And for $j=0$,
\[
\sigma_0^T\mathbf{\omega}=\sum_{t=1}^{sm+1}(m+t)t=\frac{2sm+3m+2}{3}\binom{sm+2}{2},
\]
then $\sigma_0^T\mathbf{\omega}=v_1^{(s,m)}\bar{\alpha}_1^{(s,m)}$
for $j\geq1$
\begin{align*}
\sigma_j^T\mathbf{\omega}=\sigma_j^TR_m^{(s)}X_m^{(s)}u&=\sum_{q=0}^{sm}((s+1)m+1-q)\sum_{p=0}^{sm-q}\sigma_j(p)\\
&=\sum_{p=0}^{sm}\sigma_j(p)\sum_{q=0}^{sm-p}((s+1)m+1-q)\\
&=\sum_{p=0}^{sm}\frac{(sm-p+1)((s+2)m+p+2)}{2}\sigma_j(p).
\end{align*}
Since $(p-jm)_+=0$ for $p\leq jm$,thus gives another formula
\begin{equation}
\sigma_j^T\mathbf{\omega}=\sum_{t=0}^{(s-j)m}\frac{((s-j)m-t+1)((s+2+j)m+t+2)}{2}\binom{t+j}{2j},\label{sigmaomega}
\end{equation}
we sparating the sum into two parts,
\begin{gather*}
S_1=\frac{(s+2+j)m+2}{2}\sum_{t=0}^{(s-j)m}((s-j)m-t+1)\binom{t+j}{2j},\\
S_2=\frac{1}{2}\sum_{t=0}^{(s-j)m}t((s-j)m-t+1)\binom{t+j}{2j}.
\end{gather*}
By the hockey-stick formula,
\[
S_1=\frac{(s+2+j)m+2}{2}\binom{(s-j)m+j+2}{2j+2}=\frac{(s+2+j)m+2}{2}\bar{\alpha}_{j+1}^{(s,m)}.
\]
For $S_2$,let $N=(s-j)m$ and define two series $f(x)=\sum_{t\geq0}t\binom{t+j}{2j}x^t$ and $g(x)=\sum_{k\geq0}(k+1)x^k$,then $2S_2$ is equal to the cofficient of the $N$-th of $x$ of $f(x)g(x)$.
\[
f(x)=x\frac{d}{dx}\left(\frac{x^j}{(1-x)^{2j+1}}\right)=\frac{jx^j}{(1-x)^{2j+1}}+\frac{(2j+1)x^{j+1}}{(1-x)^{2j+2}},
\]
since the identity $\sum_{t\geq0}\binom{t+j}{2j}x^t=\frac{x^j}{(1-x)^{2j+1}}$.And $g(x)=\sum_{k\geq0}(k+1)x^k=\frac{1}{(1-x)^2}$.Then
\[
f(x)g(x)=\frac{jx^j}{(1-x)^{2j+3}}+\frac{(2j+1)x^{j+1}}{(1-x)^{2j+4}}
\]
Extract the coefficients of $x^N$ of this polynomial and substitute $(s-j)m$ into $N$,we get
\[
S_2=\frac{2j(N+1)+N}{2(2j+3)}\binom{N+j+2}{2j+2}=\frac{2j((s-j)m+1)+(s-j)m}{2(2j+3)}\bar{\alpha}_{j+1}^{(s,m)}.
\]
Thus by \ref{sigmaomega},
\[
\sigma_j^TC_m^{(s)}=S_1+S_2=\frac{(j+1)((2s+3)m+3)}{2j+3}\bar{\alpha}_{j+1}^{(s,m)}.
\]
Which proved $\sigma_j^T\mathbf{\omega}=\eta_{j+1}^{(s,m)}\bar{\alpha}_{j+1}^{(s,m)}$ for $0\leq j\leq s-1$.
The second equality can be vertified parallelly.For $0\leq j\leq s-1$,let $N=(s-j)m$,
\[
\tau_j^T\mathbf{\omega}=\sum_{t=0}^{N}\frac{(N-t+1)((s+2+j)m+t+2)}{2}\binom{t+j+1}{2j+1},
\]
and suppose
\begin{equation}
\tau_j^T\mathbf{\omega}=\frac{(s+2+j)m+2}{2}\sum_{t=0}^N(N-t+1)\binom{t+j+1}{2j+1}+\frac{1}{2}\sum_{t=0}^Nt(N-t+1)\binom{t+j+1}{2j+1},\label{tauomega1}
\end{equation}
In the proof of lemma \ref{sigmatau},we have proved
\[
\sum_{t=0}^N(N-t+1)\binom{t+j+1}{2j+1}=\binom{N+j+3}{2j+3}.
\]
For the other summing part,samely defining two series $\hat{f}(x)=\sum_{t\geq0}t\binom{t+j+1}{2j+1}$ and $\hat{g}(x)=g(x)$,we need to compute the coefficient of $x^N$ of $\hat{f}(x)\hat{g}(x)$.
\[
\hat{f}(x)=x\frac{d}{dx}\left(\frac{x^j}{(1-x)^{2j+2}}\right)=\frac{jx^j}{(1-x)^{2j+2}}+\frac{(2j+2)x^{j+1}}{2j+3},
\]
then
\[
\sum_{t=0}^Nt(N-t+1)\binom{t+j+1}{2j+1}=\frac{(j+1)N+j}{j+2}\binom{N+j+3}{2j+3}.
\]
Thus by \ref{tauomega1} and substituting $(s-j)m$ into $N$,gives
\[
\tau_j^T\mathbf{\omega}=\frac{(2js+3s+3j+4)m+3j+4}{2j+4}\binom{(s-j)m+j+3}{2j+3},
\]
we have proved the identiy in lemma \ref{sigmatau},
\[
\binom{(s-j)m+j+3}{2j+3}=\bar{\alpha}_{j+1}^{(s,m)}+\bar{\beta}_{j+1}^{(s,m)}.
\]
which prove the identity about $\mathbf{\mu}$.

Next to proof the identity about $\mathbf{\nu}$.For $j=0$,
\begin{gather*}
\sigma_0^T\mathbf{\nu}=\sum_{t=1}^{sm+1}t=\binom{sm+2}{2}=\bar{\alpha}_1^{(s,m)},\\
\tau_0^T\mathbf{\nu}=\sum_{t=0}^{sm}(t+1)(sm+1-t)=\binom{sm+3}{3}=\bar{\rho}_1^{(s,m)}.
\end{gather*}
While for $1\leq j\leq s-1$,
\begin{gather*}
\sigma_j^T\mathbf{\nu}=\sum_{t=0}^{sm}(sm+1-t)\sigma_j(t)=\sum_{t=0}^{(s-j)m}((s-j)m-t+1)\binom{t+j}{2j}=\bar{\alpha}_{j+1}^{(s,m)},\\
\tau_j^T\mathbf{\nu}=\sum_{t=0}^{sm}(sm+1-t)\tau_j(t)=\sum_{t=0}^{(s-j)m}((s-j)m-t+1)\binom{t+j+1}{2j+1}=\binom{(s-j)m+j+3}{2j+3}
\end{gather*}
Since $\binom{(s-j)m+j+3}{2j+3}=\bar{\rho}_{j+1}^{(s,m)}$,we have proved lemma \ref{omega} completely.
\end{proof}

Now we proceed to verify Theorem \ref{Gmy}.
\begin{proof}[Proofs of Theorem \ref{Gmy}]
By lemma \ref{UVmsz},
\begin{gather*}
U_m^{(s)}(z)=\frac{((s+1)m+2)((s+1)m+1)}{2}+z\mathbf{\mu}^T(I-zC_m^{(s)})^{-1}\mathbf{\omega},\\
V_m^{(s)}(z)=(s+1)m+1+z\mathbf{\mu}^T(I-zC_m^{(s)})^{-1}\mathbf{\nu}.
\end{gather*}
Since for every integer $0\leq k\leq sm$,$\sigma_0(k)=1$ and $\tau_0(k)=k+1$,
\[
\sigma_0=u_{sm-1},\tau_0=(1,2,\cdots,sm+1)^T,
\]
then by $\bar{\beta}_0^{(s,m)}=(s+1)m+1$,gives
\[
\mathbf{\mu}=((s+1)m+2)\sigma_0-\tau_0=(\bar{\beta}_0^{(s,m)}+1)\sigma_0-\tau_0.
\]
by equation \ref{IABN},
\[
(I-zC_m^{(s)})^{-1}=I+C_m^{(s)}(I-zC_m^{(s)})^{-1},
\]
then gives
\begin{equation}
\mathbf{\mu}^T(I-zC_m^{(s)})^{-1}\mathbf{\omega}=(\bar{\beta}_0^{(s,m)}+1)\sigma_0^T(I-zC_m^{(s)})^{-1}\mathbf{\omega}-\tau_0^T(I-zC_m^{(s)})^{-1}\mathbf{\omega}.\label{ICsN}
\end{equation}
Let \[E_j(z)=\sigma_j^T(I-zC_m^{(s)})^{-1}\mathbf{\omega},\quad F_j=\tau_j^T(I-zC_m^{(s)})^{-1}\mathbf{\omega}.\]
We get
\[
E_j(z)=\sigma_j^T(I+zC_m^{(s)}(I-zC_m^{(s)})^{-1})\mathbf{\omega}=\sigma_j^T\mathbf{\omega}
+z\sigma_j^TC_m^{(s)}(I-zC_m^{(s)})^{-1}\mathbf{\omega},
\]
by lemma \ref{sigmatau} and identity \ref{sigmaC},
$\sigma_j^TC_m^{(s)}=\bar{\alpha}_{j+1}^{(s,m)}\sigma_0^T-\sigma_{j+1}^T$,then
\begin{align}
E_j(z)&=\sigma_j^T\mathbf{\omega}
+z(\bar{\alpha}_{j+1}^{(s,m)}\sigma_0^T(I-zC_m^{(s)})^{-1}\mathbf{\omega}-\sigma_{j+1}^T(I-zC_m^{(s)})^{-1}\mathbf{\omega})\notag\\
&=\sigma_j^T\mathbf{\omega}
+z(\bar{\alpha}_{j+1}^{(s,m)}E_0(z)-E_{j+1}(z))\label{Ejz}
\end{align}
Setting \[p_j(z):=\sum_{r=0}^{s-j-1}(-1)^r\bar{\alpha}_{j+r+1}^{(s,m)}z^r,0\leq j\leq s,\] and suppose $p_s(z)=0$.It is no hard to see that $p_j(z)+zp_{j+1}(z)=\bar{\alpha}_{j+1}^{(s,m)}$.Combining lemma \ref{omega} with \ref{Ejz},then we get
\[
E_j(z)=p_j(z)(\eta_{j+1}^{(s,m)}+zE_0(z))\quad (0\leq j\leq s),
\]
thus setting $j=0$,gives $E_0(z)=p_0(z)(\eta_1^{(s,m)}+zE_0(z))$ and
\begin{equation}
E_0(z)=\frac{p_0(z)\eta_1^{(s,m)}}{1-zp_0(z)}.\label{E0z}
\end{equation}

In the same way,we obtain
\[
F_j(z)=\tau_j(z)^T\mathbf{\omega}+z\tau_j^TC_m^{(s)}(I-zC_m^{(s)})^{-1}\mathbf{\omega},
\]
by the identity \ref{tauC},we also get
\begin{align}
F_j(z)&=\tau_j^T\mathbf{\omega}+z(\bar{\rho}_{j+1}^{(s,m)}u_{sm-1}^T(I-zC_m^{(s)})^{-1}\mathbf{\omega}-\tau_{j+1}^T(I-zC_m^{(s)})^{-1}\mathbf{\omega})\notag \\
&=\tau_j^T\mathbf{\omega}+z(\bar{\rho}_{j+1}^{(s,m)}E_0(z)-F_{j+1}(z)).\label{Fjz}
\end{align}
where $\bar{\rho}_{j+1}^{(s,m)}=\bar{\alpha}_{j+1}^{(s,m)}+\bar{\beta}_{j+1}^{(s,m)}$.
Set
\[
q_j(z):=\sum_{r=0}^{s-j-1}(-1)^r\rho_{j+1}^{(s,m)}z^r,0\leq j\leq s.
\]
where $q_s(z)=0$.Thus by lemma \ref{omega} and equation \ref{Fjz},
\[
F_j(z)+zF_{j+1}=\bar{\rho}_{j+1}^{(s,m)}(\xi_{j+1}^{(s,m)}+zE_0(z)),
\]
since $q_j(z)+zq_{j+1}(z)=\bar{\rho}_{j+1}^{(s,m)}$,then
\[
F_j(z)=q_j(z)(\xi_{j+1}^{(s,m)}+zE_0(z)).
\]
Let $j=0$ and by \ref{E0z},
\begin{equation}
F_0(z)=\frac{q_0(z)(\xi_1^{(s,m)}(1-zp_0(z))+p_0(z)\eta_1^{(s,m)})}{1-zp_0(z)}.\label{F0z}
\end{equation}
Conbine \ref{ICsN} with \ref{E0z} and \ref{F0z},
\begin{align*}
\mathbf{\mu}^T(I-zC_m^{(s)})^{-1}\mathbf{\omega}&=(\bar{\beta}_0^{(s,m)}+1)E_0(z)-F_0(z)\\&=\frac{(\bar{\beta}_0^{(s,m)}+1-q_0(z))p_0(z)\eta_1^{(s,m)}-q_0(z)\xi_1^{(s,m)}(1-zp_0(z))}{1-zp_0(z)}
\end{align*}
Then
\[
U_m^{(s)}(z)=\frac{\bar{\beta}_0^{(s,m)}(\bar{\beta}_0^{(s,m)}+1)}{2}+\frac{z((\bar{\beta}_0^{(s,m)}+1-q_0(z))p_0(z)\eta_1^{(s,m)}-q_0(z)\xi_1^{(s,m)}(1-zp_0(z)))}{1-zp_0(z)}.
\]

The formula of $V_m^{(s)}(z)$ can be proved parallelly.By \ref{IABN},
\begin{equation}
\mathbf{\mu}^T(I-zC_m^{(s)})^{-1}\mathbf{\nu}=(\bar{\beta}_0^{(s,m)}+1)\sigma_0^T(I-zC_m^{(s)})^{-1}\mathbf{\nu}-\tau_0^T(I-zC_m^{(s)})^{-1}\mathbf{\nu},\label{muCnu}
\end{equation}
Samely let $\bar{E}_j(z)=\sigma_j^T(I-zC_m^{(s)})^{-1}\mathbf{\nu}$ and $\bar{F}_j(z)=\tau_j^T(I-zC_m^{(s)})^{-1}\mathbf{\nu}$.Decomposing $(I-zC_m^{(s)})^{-1}$ and applying lemma \ref{sigmatau} gives
\begin{gather*}
\bar{E}_j(z)=\sigma_j^T\mathbf{\nu}
+z(\bar{\alpha}_{j+1}^{(s,m)}u_{sm-1}^T-\sigma_{j+1}^T)(I-zC_m^{(s)})^{-1}\mathbf{\nu}=\bar{\alpha}_{j+1}^{(s,m)}
+z(\bar{\alpha}_{j+1}^{(s,m)}\bar{E}_0(z)-\bar{E}_{j+1}(z)),\\
\bar{F}_j(z)=\tau_j^T\mathbf{\nu}
+z(\bar{\rho}_{j+1}^{(s,m)}u_{sm-1}^T-\tau_{j+1}^T)(I-zC_m^{(s)})^{-1}\mathbf{\nu}=\bar{\rho}_{j+1}^{(s,m)}
+z(\bar{\rho}_{j+1}^{(s,m)}\bar{F}_0(z)-\bar{F}_{j+1}(z)),
\end{gather*}
then we get
\[
\bar{E}_j(z)+z\bar{E}_{j+1}(z)=\bar{\alpha}_{j+1}^{(s,m)}(1+z\bar{E}_0(z)),\quad
\bar{F}_j(z)+z\bar{F}_{j+1}(z)=\bar{\rho}_{j+1}^{(s,m)}(1+z\bar{E}_0(z)).
\]
Therefore by the recurrence of the polynomials $p_j(z)$ and $q_j(z)$,
\[
\bar{E}_j(z)=p_j(z)(1+z\bar{E}_0(z)),\quad \bar{F}_j(z)=q_j(z)(1+z\bar{E}_0(z)).
\]
Setting $j=0$,gives
\[
\bar{E}_0(z)=\frac{p_0(z)}{1-zp_0(z)},\quad \bar{F}_0(z)=\frac{q_0(z)}{1-zp_0(z)}.
\]
Then
\[
V_m^{(s)}(z)=\bar{\beta}_0^{(s,m)}+\frac{z((\bar{\beta}_0^{(s,m)}+1)p_0(z)-q_0(z))}{1-zp_0(z)}.
\]
Recalls that the definition of $G_m^{s}(y)$,and separate into odd- and even- parts,
\[
G_m^{(s)}(y)=\sum_{r\geq0}L_{2r}^{(s)}(m)y^{2r}+\sum_{r\geq0}L_{2r+1}^{(s)}(m)y^{2r+1}=1+y^2U_m^{(s)}(y^2)+yV_m^{(s)}(y^2).
\]
This the first-half complete proof of theorem \ref{Gmy}
\end{proof}
\begin{flushleft}
\textbf{Remark.}\emph{For a special case:$m=1$,then $L_d^{(s)}(1)=\sharp(\mathcal{P}^{(s)}\cap \mathbb{Z}^d)=N_d(s)$.By the definition of $G_m^{(s)}(y)$ and $F(x)$,
\begin{align*}
G_1^{(s)}(y)=\sum_{d\geq0}L_d^{(s)}(1)y^d&=L_0^{(s)}(1)+L_1^{(s)}(1)y+\sum_{d\geq2}L_d^{(s)}(1)y^d\\
&=1+(s+2)y+F(y).
\end{align*}}
\end{flushleft}
\subsection{The Explicit Expression for the Exceptional Case $s=1$}
The case \(s=1\) is exceptional, as the dimension generating function \(G_m^{(1)}(y)\) admits a rational form. For \(m\geq 0\), recall that
\[
G_m^{(1)}(y)=\sum_{d\geq0}L_d^{(1)}(m)y^d,L_d^{(1)}(m)=\sharp(m\mathcal{P}_d^{(1)}\cap \mathbb{Z}^d).
\]
We utilize three auxiliary sets according to the constraints of these inequality system,let
\begin{gather*}
\Delta_m=\{(x,y)\in \mathbb{R}_{\geq0}^2:x+y\leq m\},\\
\Gamma_m=\{(x,y,z)\in \mathbb{R}_{\geq0}^3:x+y\leq 2m,y+z\leq m\},\\
\Theta_m=\{(x,y,z)\in \mathbb{R}_{\geq0}^3:x+y\leq m,y+z\leq2m\}.
\end{gather*}
Write the number of integer points of these sets
\[
\delta_m=\sharp(\Delta_m \cap \mathbb{Z}^2),\gamma_m=\sharp(\Gamma_m \cap \mathbb{Z}^3),\theta_m=\sharp(\Theta_m \cap \mathbb{Z}^3).
\]
We find that these formula have the follows form,
\begin{gather}
\delta_m=\sum_{x=0}^m(m+1-x)=\binom{m+2}{2},\label{delta}\\ 
\gamma_m=\theta_m=\sum_{y=0}^m(m+1-y)(2m+1-y)=\frac{(m+1)(m+2)(5m+3)}{6}.\label{gamma}
\end{gather}
Preserving the geometric structure of the constraints defining \(m\mathcal{P}_d^{(1)}\), consider the even-dimensional case \(d=2r\), where \(r\geq 1\). For \(1\leq j\leq r\), each odd-indexed constraint is of the form $x_{2j-1}+x_{2j}\leq 2m$,whereas each even-indexed constraint is of the form $x_{2j}+x_{2j+1}\leq m$.In particular, the latter implies $x_{2j},x_{2j+1}\leq m$ for \(j\geq 1\). Consequently, for \(1\leq j\leq r-2\), the constraint $x_{2j+1}+x_{2j+2}\leq 2m$ is automatically satisfied, since
$$
x_{2j+1}\leq m
\qquad\text{and}\qquad
x_{2j+2}\leq m
$$
follow from the constraints
$$
x_{2j}+x_{2j+1}\leq m
\qquad\text{and}\qquad
x_{2j+2}+x_{2j+3}\leq m,
$$
respectively. Thus, among the original \(2r-1\) constraints, the constraints of the form \(x_{2j+1}+x_{2j+2}\leq 2m\) for \(1\leq j\leq r-2\) are redundant. We can therefore reduce the \(2r-1\) constraints to the following equivalent system,
\[
x_1+x_2\leq2m,x_2+x_3\leq m,x_{2j}+x_{2j+1}\leq m,x_{2r-2}+x_{2r-1}\leq m,x_{2r-1}+x_{2r}\leq2m,
\]
for $2\leq j\leq r-2$.

Then exsits a bijection $\pi$ such that
\[
m\mathcal{P}_{2r}^{(1)}\cong \Gamma_m\times \Delta_m^{r-3}\times \Theta_m\quad(r\geq3),
\]
Consequently,by formulas \ref{delta} and \ref{gamma}
\[
L_{2r}^{(1)}(m)=\sharp(m\mathcal{P}_{2r}^{(1)}\cap \mathbb{Z}^{2r})=\delta_m^{r-3}\gamma_m^2  \quad (r\geq3).
\]
For the odd-dimensional case \(d=2r+1\), we obtain the corresponding system by adding the even-indexed constraint $x_{2r}+x_{2r+1}\leq m$ to the system considered in the even-dimensional case. In this setting, the last odd-indexed constraint $x_{2r-1}+x_{2r}\leq 2m$ is automatically satisfied whenever both $x_{2r-2}+x_{2r-1}\leq m$ and $x_{2r}+x_{2r+1}\leq m$ hold. Indeed, these two constraints imply $x_{2r-1}\leq m$ and $x_{2r}\leq m$,and hence $x_{2r-1}+x_{2r}\leq 2m$.Thus, the last odd-indexed constraint is redundant, and the original system can be simplified in the same manner as in the even-dimensional case. Consequently, the constraints reduce to
\begin{align*}
x_1+x_2\leq2m,x_2+x_3\leq m.x_{2j}+x_{2j+1}\leq m,2\leq j\leq r,
\end{align*}
and exists a bijection such that
\begin{align*}
m\mathcal{P}_{2r+1}^{(1)}\cong \Gamma_m\times \Delta_m^{r-1} \quad (r\geq1).
\end{align*}
therefore
\begin{align*}
L_{2r+1}^{(1)}(m)=\sharp(m\mathcal{P}_{2r+1}^{(1)}\cap \mathbb{Z}^{2r})=\delta_m^{r-1}\gamma_m \quad (r\geq1)
\end{align*}
Since the definitons of $L_0^{(s)}(m)=1$ and $L_1^{(s)}(m)=(s+1)m+1$,
\begin{align}
G_m^{(1)}(y)&=\sum_{d\geq0}L_d^{(1)}(m)y^d=\sum_{r\geq0}L_{2r}^{(1)}(m)y^{2r}+\sum_{r\geq0}L_{2r+1}^{(1)}(m)y^{2r+1}\notag \\
&=L_0^{(1)}(m)+L_1^{(1)}(m)y+L_2^{(1)}(m)y^2+L_4^{(1)}(m)y^4+\sum_{r\geq3}\delta_m^{r-3}\gamma_m^2y^{2r}
+\sum_{r\geq1}\delta_m^{r-1}\gamma_my^{2r+1}.\label{Gm1y}
\end{align}
where $L_0^{(1)}(m)=1,L_1^{(1)}(m)=2m+1$ and
\begin{gather*}
L_2^{(1)}(m)=\sharp\{(x_1,x_2)\in \mathbb{Z}^2:x_1+x_2\leq2m\}=\binom{2m+2}{2},\\
L_4^{(1)}(m)=\sum_{x_2+x_3\leq m}(2m+1-x_2)(2m+1-x_3)=\frac{(m+1)(m+2)(11m^2+13m+4)}{8}.
\end{gather*}
The second identity since the nonnagetive positive integer solution  of system $x_1+x_2\leq2m,x_2+x_3\leq m,x_3+x_4\leq2m$.Then simplify \ref{Gm1y} gives
\[
G_m^{(1)}(y)=\frac{1+(2m+1)y+3\binom{m+1}{2}y^2-\binom{m+2}{3}y^3+\frac{9m+3}{4}\binom{m+2}{3}y^4-\frac{1}{3}\binom{m+2}{2}\binom{m+3}{4}y^6}{1-\binom{m+2}{2}y^2}.
\] 
Thus,the complete proof of this theorem \ref{Gmy} is now finished.
\section{The $h^*$-series and Ehrhart series}
By lemma \ref{Hksy},we have the first four series of $H_k^{(s)}(y)$ for fixing $s=1$,
\begin{align*}
H_0^{(1)}(y)&=G_0^{(1)}(y)=\frac{1}{1-y},\\
H_1^{(1)}(y)&=G_1^{(1)}(y)-\frac{d}{dy}(yG_0^{(1)}(y))=\frac{y(y^7-2y^6-2y^5+7y^4-8y^3+4y^2-y-1)}{(3y^2-1)(1-y)^2},\\
H_2^{(1)}(y)&=G_2^{(1)}(y)-\frac{d}{dy}(yG_1^{(1)}(y))+\frac{y}{(1-y)^3}\\
&=\frac{y^4(30y^7-90y^6+105y^5-36y^4+24y^3-18y^2+21y-10)}{(6y^2-1)(3y^2-1)(1-y)^3},\\
H_3^{(1)}(y)&=G_3^{(1)}(y)-\frac{d}{dy}(yG_2^{(1)}(y))+\frac{y}{2}\frac{d^2}{dy^2}(yG_1^{(1)}(y))-\frac{y^2}{(1-y)^4}\\&=\frac{-y^5R_3^{(1)}(y)}{(10y^2-1)(6y^2-1)^2(3y^2-1)^3(1-y)^4},\\
H_4^{(1)}(y)&=G_4^{(1)}(y)-\frac{d}{dy}(yG_3^{(1)}(y))+\frac{y}{2}\frac{d^2}{dy^2}(yG_2^{(1)})-\frac{y^2}{6}\frac{d^3}{dy^3}(yG_1^{(1)}(y))+\frac{y^3}{(1-y)^5}\\
&=\frac{y^6R_4^{(1)}(y)}{(15y^2-1)(10y^2-1)^2(6y^2-1)^3(3y^2-1)^4(-1+y)^5}.
\end{align*}
where
\begin{align*}
R_3^{(1)}(y)&=4500y^{11}-12168y^{10}+22704y^9-17052y^8+9109y^7-8752y^6+2371y^5\\&+1328y^4-405y^3+328y^2-159y-4,\\
R_4^{(1)}(y)&=8262000y^{19}-27896400y^{18}+73677600y^{17}-91436040y^{16}+
84592080y^{15}-85465476y^{14}\\&+37676340y^{13}+14097186y^{12}-13827099y^{11}
+6845550y^{10}-5455411y^9-43612y^8\\&+1568426y^7-289840y^6-29722y^5-
2382y^4-8095y^3+2838y^2+41y+16.
\end{align*}
and by identity $\sum_{d\geq0}Ehr_{\mathcal{P}_d^{(s)}}(z)y^d=\sum_{m\geq0}G_m^{(s)}(y)z^m$,we have the concrete forms of the Ehrhart series,
\begin{align*}
Ehr_{\mathcal{P}_0^{(1)}}(z)&=\frac{1}{1-z},& Ehr_{\mathcal{P}_1^{(1)}}(z)&=\frac{1+z}{(1-z)^2},\\
 Ehr_{\mathcal{P}_2^{(1)}}(z)&=\frac{1+3z}{(1-z)^3},&
Ehr_{\mathcal{P}_3^{(1)}}(z)&=\frac{1+4z}{(1-z)^4},\\ Ehr_{\mathcal{P}_4^{(1)}}(z)&=\frac{1+16z+16z^2}{(1-z)^5},& Ehr_{\mathcal{P}_5^{(1)}}(z)&=\frac{1+18z+27z^2+4z^3}{(1-z)^6}.
\end{align*}
\begin{flushleft}
\textbf{Remark:}\emph{Combining Ehrhart theory with the small-dimensional cases analyzed above, we conclude that the $h^*$-coefficients are non-negative and the $h^*$-polynomials are not Gorenstein across all dimensions. However, unimodality may still hold for a fixed $s=1$. Selected values of $h_{d,k}^{(1)}$ are presented in Table\ref{tab}.}
\end{flushleft}
\begin{table}[H]
\centering
\label{tab}
    \begin{tabular}{cc|ccccccc}
        \toprule[1.5pt]
        \multirow{2}{*}{$\mathbf{h_{d,k}^{(1)}}$} & & \multicolumn{7}{c}{$\mathbf{k}$} \\
        \cmidrule{3-9}
        & & 0 & 1 & 2 & 3 & 4 & 5 & 6\\
        \midrule[1pt]
        \multirow{10}{*}{$\mathbf{d}$} 
        & 0 & 1 & 0 & 0 & 0 & 0 & 0 & 0\\
    
        & 1 & 1 & 1 & 0 & 0 & 0 & 0 & 0\\

        & 2 & 1 & 3 & 0 & 0 & 0 & 0 & 0\\

        & 3 & 1 & 4 & 0 & 0 &0 & 0 & 0\\

        & 4 & 1 & 16 & 16 & 0 & 0 & 0 &0 \\

        & 5 & 1 & 18 & 27 & 4 &0 & 0 & 0\\

        & 6 & 1 & 55 & 243 & 175 &16 & 0 & 0\\

        & 7 & 1 & 64 & 388 & 472 & 121 & 4 & 0\\
 
        & 8 & 1 & 177 & 2298 & 61888 &4323 & 717 & 16\\
        \bottomrule[1.5pt]
    \end{tabular}
\caption{some first value of $h_{d,k}^{(1)}$}
\end{table}
\textbf{Acknowledgments:}The auther would like to thank Professor Yueming Zhong for functional suggestions and a M.S's encouragement.

\end{document}